\documentclass{article}
\usepackage{PRIMEarxiv}
\usepackage[colorlinks=true,linkcolor=red,urlcolor=blue,citecolor=blue,bookmarks=true]{hyperref}
\usepackage{lipsum}
\usepackage{amsfonts}
\usepackage{subfig}
\usepackage{graphicx}
\usepackage{epstopdf}
\usepackage{amsmath}
\usepackage{amssymb}
\usepackage{cleveref}
\usepackage{amsopn}
\usepackage{setspace}
\usepackage{color}
\usepackage{algorithmicx}
\usepackage[ruled]{algorithm}
\usepackage{algpseudocode}
\usepackage{algpascal}
\usepackage{listings}
\usepackage{footnotebackref}
\usepackage{threeparttable, tablefootnote}
\usepackage[table]{xcolor}
\usepackage{amsthm}
\usepackage{latexsym}
\usepackage[mathscr]{euscript}
\def\bxi{\mbox{\boldmath $\xi$}}

\DeclareMathAlphabet\mathbfcal{OMS}{cmsy}{b}{n}
\usepackage[utf8]{inputenc} 
\usepackage[T1]{fontenc}    
\usepackage{url}            
\usepackage{booktabs}        
\usepackage{nicefrac}       
\usepackage{microtype}      
\usepackage{fancyhdr} 

\lstdefinestyle{customc}{
  belowcaptionskip=1\baselineskip,
  breaklines=true,
  frame=single,
  xleftmargin=\parindent,
  language=Fortran,
  showstringspaces=false,
  basicstyle=\footnotesize\ttfamily,
  keywordstyle=\bfseries\color{green!40!black},
  commentstyle=\itshape\color{purple!40!black},
  identifierstyle=\color{blue},
  stringstyle=\color{orange},
}

\lstdefinestyle{customm}{
  belowcaptionskip=1\baselineskip,
  breaklines=true,
  frame=single,
  xleftmargin=\parindent,
  language=Matlab,
  showstringspaces=false,
  basicstyle=\footnotesize\ttfamily,
  keywordstyle=\bfseries\color{green!40!black},
  commentstyle=\itshape\color{purple!40!black},
  identifierstyle=\color{blue},
  stringstyle=\color{orange},
}

\lstdefinestyle{customp}{
  belowcaptionskip=1\baselineskip,
  breaklines=true,
  frame=single,
  xleftmargin=\parindent,
  language=Python,
  showstringspaces=false,
  basicstyle=\footnotesize\ttfamily,
  keywordstyle=\bfseries\color{green!40!black},
  commentstyle=\itshape\color{purple!40!black},
  identifierstyle=\color{blue},
  stringstyle=\color{orange},
}

\graphicspath{{media/}}     

\title{Solving Stochastic PDEs using FEniCS and UQTk}

\author{
  Ajit Desai \\
  Department of Civil and Environmental Engineering \\
  Carleton University, Ottawa, ON, Canada \\
  \texttt{ajit.ndesai@gmail.com} \\
}

\begin{document}
\maketitle

\begin{abstract}
The intrusive (sample-free) spectral stochastic finite element method (SSFEM) is a powerful numerical tool for solving stochastic partial differential equations (PDEs). However, it is not widely adopted in academic and industrial applications because it demands intrusive adjustments in the PDE solver, which require substantial coding efforts compared to the non-intrusive (sampling) SSFEM. 
Using an example of stochastic PDE, in this article, we demonstrate that the implementational challenges of the intrusive approach can be alleviated using FEniCS---a general purpose finite element package and UQTk---a  collection of libraries and tools for the quantification of uncertainty. Furthermore, the algorithmic details and code snippets are provided to assist computational scientists in implementing these methods for their applications. This article is extracted from the author's thesis~\cite{desai2019scalable}.
\end{abstract}


\section{Introduction}\label{sec:intro}
This article is extracted from the appendices of the author's Ph.D. thesis~\cite{desai2019scalable}.
The spectral stochastic finite element method (SSFEM) is a powerful numerical tool employed for uncertainty quantification (UQ) of stochastic partial differential equations (PDEs)~\cite{ghanemSFEM1991,le2010spectral}.
The SSFEM is based on polynomial chaos expansion (PCE), i.e., a series representation of random vectors in terms of orthogonal polynomials~\cite{ghanemSFEM1991,le2010spectral}.

SSFEM is developed by leveraging the advantages of the deterministic finite element method (FEM), and its application requires the following three steps:
(1) spatial discretization of a stochastic PDE using a FEM, (2) stochastic discretization of the random system parameters, stochastic source term, and solution process using the PCE, followed by a standard Galerkin projection along the random dimensions, and (3) the resulting system is solved for the PCE coefficients of the solution process using an intrusive (non-sampling) or non-intrusive (sampling) approach~\cite{ghanemSFEM1991,le2010spectral,ghanem1999ingredients}.

Using an example of stochastic PDE, this article largely focuses on the implementation aspects of the above three steps. In step 1, we use the FEniCS---a general purpose finite element package for spatial discretization of a PDE using a FEM~\cite{fenics/dolfin17}. 
In step 2, we use the UQTk---a collection of libraries and tools for uncertainty quantification~\cite{debusschere2013uqtk}, for the stochastic discretization of the random system parameters solution process using the PCE. Finally, in step 3, we use FEniCS provided solver to solve both: (a) a coupled system of equations arising in the context of an intrusive approach and (b) an individual sample for a non-intrusive approach.

There are numerous articles present in the literature that provide in-depth details for the formulation and solution of both intrusive and non-intrusive SSFEM~\cite{ghanemSFEM1991,le2010spectral,ghanem1999ingredients,eldred2009comparison,xiu2003modeling,desai2010analysis}. Moreover, many researchers are focused on developing domain decomposition-based algorithms in conjunction with high-performance computing to efficiently tackle stochastic PDEs using SSFEM~\cite{sarkarIJNME2009,subber2012PhDTh,elman2007solving,mandel2008multispace,ghosh2009feti,stavroulakis2017gpu,pellissetti2000iterative,sousedik2014hierarchical,papadrakakis2011new,desai2017scalable,desai2022domain}. However, this article focuses on solving a stochastic PDE using serial solver in Python on a traditional desktop computer.

This remaining article is organized in the following manner. \Cref{sec:klepce} is dedicated to the spectral representation of input and output stochastic processes using the polynomial chaos expansion (PCE). This is followed by the formulation and implementation of the intrusive and non-intrusive SSFEM in the \Cref{sec:ssfem}. Finally, in \Cref{sec:conclusion} we conclude our findings. Various code snippets are provided in~\Cref{sec:sinppets}, to assist computational scientists in implementing these methods for their applications. For additional details refer to~\cite{desai2019scalable}.


\section{Spectral Representation of Stochastic Process}\label{sec:klepce}

The most widely utilized approaches for the spectral representation of input and output stochastic processes are the Karhunen-Lo{\`e}ve expansion (KLE) and polynomial chaos expansion (PCE), which we briefly discuss next. For further details refer to~\cite{desai2019scalable} and many articles cited therein.

\subsection{Karhunen-Lo{\`e}ve Expansion} \label{apn:kle}

Consider $\alpha (\textbf{\textit{x}},  {\bxi}(\theta))$ to be a real-valued stochastic process, a function of the position vector $\textbf{\textit{x}}$ defined over physical domain $\mathcal{D}$ and the set of random variables $\bxi$ which are a function of random event $\theta$ defined by complete probability space ($\Omega, \mathcal{E}, \mathcal{P}$). The KL expansion of an arbitrary non-Gaussian and non-stationary stochastic process using $L$ random variables can be written as~\cite{ghanemSFEM1991,ghanem1999ingredients},
\begin{equation}\label{eq:apnKLE}
\alpha (\textbf{\textit{x}},  {\bxi}(\theta)) = \bar{\alpha}(\textbf{\textit{x}}) + \sum^{L}_{n=1} \sqrt{{\lambda}_n} \ f_n(\textbf{\textit{x}}) {\xi_{n}}(\theta),
\end{equation}
where $\bar{\alpha}(\textbf{\textit{x}}) $ is the expected value of the random process, \{${\xi_{n}}$\} is a set of uncorrelated (not necessarily independent) random variables, \{${\lambda}_n$\} and \{$f_n$\} are the eigenvalues and eigenfunctions of the covariance function $C_{\alpha\alpha}({\textbf{\textit{x}}}, {\textbf{\textit{y}}})$ obtained by solving the following integral equation~\cite{ghanemSFEM1991}
\begin{equation}
\int_{\mathcal{D}} C_{\alpha \alpha}( \textbf{\textit{x}}, \textbf{\textit{y}})f_n(\textbf{\textit{y}})d\textbf{\textit{y}} = \lambda_{n} f_n(\textbf{\textit{x}}).
\label{eq:exp2d_int}
\end{equation}

For example, consider an exponential covariance function of a stochastic process defined over a square domain $\mathcal{D}(x,y)$
over the interval $[-a \ \ a] \times [-a \ \ a]$~\cite{ghanemSFEM1991}, 
\begin{equation}
C({x}_1, {y}_1; {x}_2, {y}_2) = {\sigma}^2 e^{- | x_2 - x_1 | / b_x \ - | y_2 - y_1 | / b_y},
\label{eq:apn_exp2d}
\end{equation}
using $b_x = b_y = b = 1$, the correlation lengths along $x$ and $y$ directions respectively and $\sigma^2$ denotes the variance of the stochastic process.

Solving the integral equation given in \Cref{eq:exp2d_int} for the covariance kernel in \Cref{eq:apn_exp2d}, the eigenvalues and eigenfunctions are obtained as~\cite{ghanemSFEM1991},
\begin{equation}\label{eq:apn_eigval_2D}
  \lambda_n = \lambda^{x}_i  \otimes \lambda^{y}_i ,
\end{equation}
  \begin{equation}
f_n(x, y) = g_i (x) \otimes h_i(y).
\label{eq:apn_eigfn_2D}
\end{equation}
where $\otimes$ denotes the tensor product (for example, see~\Cref{table:eigenTensor}). For simplicity, we refer $\{ \lambda^{x}_i, \lambda^{y}_i \}$ and $\{ g_i (x), h_i(y) \}$ are one-dimensional eigenvalues and eigenvectors.
\begin{align}\label{eq:apn_eigenVal}
\lambda^x_i = \lambda^y_i  &= {\sigma}  \frac{2b}{1 + b^2 {\omega_i}^2},
\end{align}
and
\begin{equation}
 g_i(\textbf{\textit{z}}) = h_i(\textbf{\textit{z}}) =
\begin{cases}
& \frac{\cos({\omega}_i\textbf{\textit{z}})}{\sqrt{a + \frac{\sin(2{\omega}_i a)}{2 {\omega}_i}}}, \ \  \text{for} \  i \ \text{odd},   \\
 & \frac{\sin({\omega}_i\textbf{\textit{z}})}{\sqrt{a - \frac{\sin(2{\omega}_i a)}{2 {\omega}_i}}}, \ \  \text{for} \  i \ \text{even}.
 \end{cases}
\label{eq:eigenVecApn}
\end{equation}
Here $\omega_i$'s are the solution of the following transcendental equations~\cite{ghanemSFEM1991},
\begin{align}
 \frac{1}{b} - {\omega_i} \ \tan(\omega_i \  a) &= 0,  \ \  \text{for} \  i \ \text{odd},  \nonumber \\
 {\omega_i} + \frac{1}{b} \ \tan(\omega_i \  a) &= 0, \ \  \text{for} \  i \ \text{even}.
 \label{eq:apn_2D_omegas}
\end{align}


\rowcolors{2}{gray!10}{white}
\begin{table}[htbp]
\onehalfspacing
\begin{center}
\begin{threeparttable}
\caption{Tensor product of one-dimensional eigenvalues.}
\begin{tabular}{|c|*{2}{c|}r}
\hline
two-dimensional eigenvalues	& one-dimensional eigenvalues  \\ \hline
$\lambda_1$	& $\lambda^{x}_1 \times \lambda^{y}_1$  \\ \hline
$\lambda_2$ & $\lambda^{x}_1 \times \lambda^{y}_2$ \\ \hline
$\lambda_3$ & $\lambda^{x}_2 \times \lambda^{y}_1$ \\
\hline
$\lambda_4$ & $\lambda^{x}_2 \times \lambda^{y}_2$ \\ \hline
\end{tabular}
\label{table:eigenTensor}
\end{threeparttable}
\end{center}
\end{table}

\begin{table}[htbp]
\onehalfspacing
\begin{center}
\begin{threeparttable}
\caption{One-dimensional omegas $\omega_i$ and eigenvalues $\lambda_i^x, i=1,2,\dots,7$.}
\begin{tabular}{|l|*{7}{c|}r}
\hline
$index$ $i$ & 1  &  2  &  3  & 4  &  5 &  6 & 7  \\ \hline
$\omega_i$	& 1.306  &  3.673  &  6.585  & 12.723  &  15.834 &  18.955 & 22.082  \\ \hline
$\lambda_i^x$ & 0.7388  &  0.1380  &  0.0451  &  0.0213  &  0.0123  &  0.0079  &  0.0056   \\
\hline
\end{tabular}
\label{table:eigen1D}
\end{threeparttable}
\end{center}
\end{table}

\begin{table}[htbp!]
\onehalfspacing
\begin{center}
\begin{threeparttable}
\caption{Two-dimensional eigenvalues $\lambda_n, n=1,2,\dots,7$ and $sortIndex_i = \{ 1,1,2,1,3,2,1 \}$ and $sortIndex_j = \{ 1,2,1,3,1,2,4 \}$.}
\begin{tabular}{|l|*{7}{c|}r}
\hline
$index$ $n$	& 1  &  2  &  3  & 4  &  5 &  6 & 7  \\ \hline
$sortIndex_i$	&  $\lambda_1^x$  &   $\lambda_1^x$   &  $\lambda_2^x$   &  $\lambda_1^x$   &  $\lambda_3^x$  &   $\lambda_2^x$   &  $\lambda_1^x$  \\ \hline
$sortIndex_j$	&  $\lambda_1^y$  &   $\lambda_2^y$   &  $\lambda_1^y$   &  $\lambda_3^y$   &  $\lambda_1^y$  &   $\lambda_2^y$   &  $\lambda_4^y$   \\ \hline
$\lambda_n$   & 0.5458  &  0.1020  &  0.1020 &   0.0333  &  0.0333    & 0.0190  &  0.0158   \\
\hline
\end{tabular}
\label{table:eigen2D}
\end{threeparttable}
\end{center}
\end{table}

\begin{figure}[htbp]
 \centering
 \subfloat[one-dimensional eigenvalues (\Cref{eq:apn_eigenVal}) \label{subfig:ev1D}]{%
 \includegraphics[width=0.48\textwidth,height=0.28\textheight]{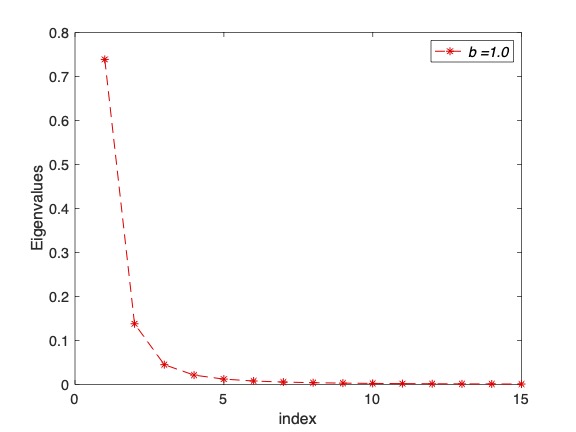}}
 \subfloat[two-dimensional eigenvalues (\Cref{eq:apn_eigval_2D})  \label{subfig:ev2D}]{%
 \includegraphics[width=0.48\textwidth,height=0.28\textheight]{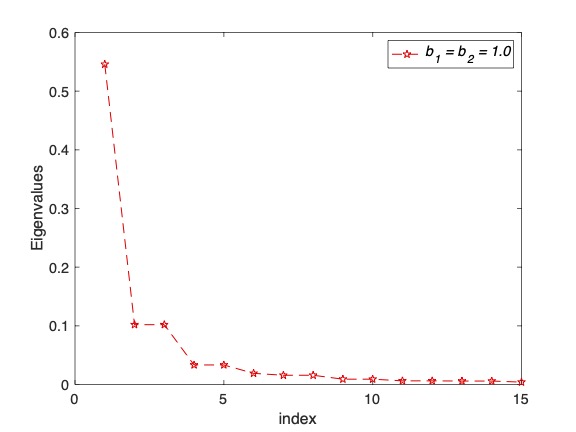}}
 \caption{Eigenvalues for $b=1.0$ and $a=0.5$}
 \label{fig:eigevalues}
\end{figure}

Solving for one-dimensional $\omega_i$ and $\{ \lambda^{x}_i, \lambda^{y}_i \}$ from \Cref{eq:apn_2D_omegas} and \Cref{eq:apn_eigenVal} respectively for $b=1$ and $a=0.5$, i.e., using unit square domain, we get the results which are summarized  in \Cref{table:eigen1D}.
Note that these results from eigenvalue analysis are obtained by sorting eigenvalues $\{\lambda_i\}_{i=1}^{7}$ in descending order (i.e., largest to smallest eigenvalues).
The first few eigenvalues account for most of the contribution to the variance and the contribution of higher indexed eigenvalues decreases quickly as shown in~\Cref{fig:eigevalues}.
The two-dimensional eigenvalues are obtained by taking tensor product of (sorted) one-dimensional eigenvalues as shown in~\Cref{table:eigenTensor}. After taking the tensor product the two-dimensional eigenvalues are sorted in descending order are shown in \Cref{table:eigen2D}. This leads us to the new one-dimensional $index$ in $x$ and $y$ dimensions, which we call
$sortIndex_i$ and $sortIndex_j$ as shown in~\Cref{table:eigen2D}.

\begin{figure}[htbp]
 \centering
 \subfloat[one-dimensional case \label{subfig:ev1D_2}]{%
 \includegraphics[width=0.48\textwidth,height=0.28\textheight]{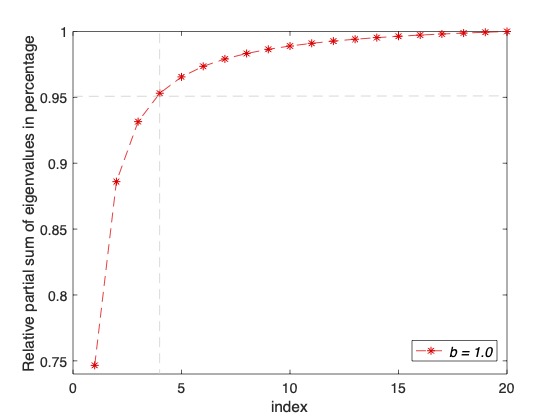}}
 \subfloat[two-dimensional case \label{subfig:ev2D_2}]{%
 \includegraphics[width=0.48\textwidth,height=0.28\textheight]{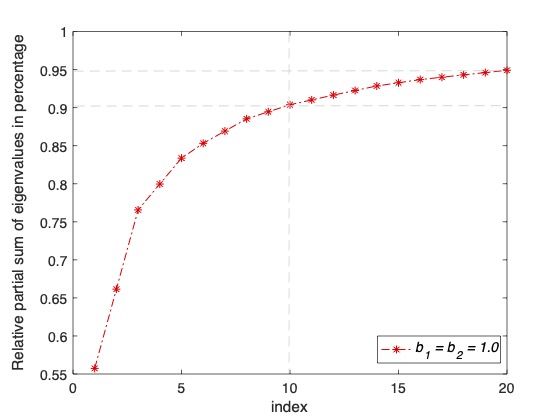}}
 \caption{Relative partial sum of eigenvalues for $b=1.0$ and $a=0.5$}
 \label{fig:eigenContribution}
\end{figure}

From \Cref{fig:eigevalues} and \ref{fig:eigenContribution}, it can be observed that the eigenvalue contribution decreases more rapidly in two-dimensional cases compared to one-dimensional cases. For example, to get the relative partial sum of eigenvalues $\Big(\frac{\sum_i^k \lambda_i}{\sum_i^n \lambda_i}, n >> k \Big)$ of $95\%$, we need $20$ eigenvalues in the two-dimensional case (\Cref{eq:apn_eigval_2D}) as oppose to only $4$ modes in the one-dimensional case (\Cref{eq:apn_eigenVal}). Therefore, the number of random variables required to characterize the underlying stochastic process can increase with the physical dimension of the problem~\cite{ghanemSFEM1991,ghanem1999ingredients}.

\subsection{Polynomial Chaos Expansion} \label{apn:pce}
Consider a random process $\alpha (\textbf{\textit{x}}, \bxi (\theta) )$, as
function of position vector $\textbf{\textit{x}}$ and set of random variables $\bxi$ which are function of a random event $\theta$. Using polynomial chaos expansion (PCE)
the stochastic process can be written as~\cite{ghanemSFEM1991} (for notational convenience $\theta$ is dropped from henceforth),
\begin{align} \label{eq:polc}
\alpha (\textbf{\textit{x}},  {\bxi} )   & =  \widehat{\alpha}_0 \Phi_0  +  \sum_{i_1 = 1}^{\infty} \widehat{\alpha}_{i_1} \Phi_ 1 (  \xi_{i_1}  ) \\ \nonumber & + \ \sum_{i_1 = 1}^{\infty} \sum_{i_2 = 1}^{i_1} \widehat{\alpha}_{i_1 i_2} \Phi_2 ( \xi_{i_1} , \xi_{i_2} )
\\ \nonumber & + \ \sum_{i_1 = 1}^{\infty} \sum_{i_2 = 1}^{i_1} \sum_{i_3 = 1}^{i_2} \widehat{\alpha}_{i_1 i_2 i_3} \Phi_3 ( \xi_{i_1}  , \xi_{i_2} , \xi_{i_3} )
 + \ldots \ ,
\end{align}
where the $\Phi_n (  \xi_{i_1} , \ldots , \xi_{i_n} ) $ are the multidimensional polynomial basis or polynomial chaoses~\cite{ghanemSFEM1991} of order $n$
in terms of $n$-dimensional random variables $( \xi_{i_1} , \ldots , \xi_{i_n} ) $. 
In \Cref{eq:polc}, $\{ \widehat{\alpha} \}$ are the deterministic PC coefficients which are function of $\textbf{\textit{x}}$.

For numerical implementation, a concise and truncated PC expansion is used. For instance, using $P_u$ terms, the PCE of $\alpha (\textbf{\textit{x}}, \bxi)$ can be written as~\cite{ghanemSFEM1991},
\begin{equation}\label{eq:pce_rewrite}
\alpha (\textbf{\textit{x}},  {\bxi})  \approx \sum_{j=0}^{P_u} \alpha_j (\textbf{\textit{x}}) \Psi_j (  \bxi ).
\end{equation}
There is an one-to-one relationship between the $\{ \Phi_i \}$ and
$\{\Psi_j \}$ and also $\{\widehat{\alpha}_i\}$ and $\{\alpha_j\}$ in
\Cref{eq:polc} and \Cref{eq:pce_rewrite}. Note that, in this analysis, $\bxi$ is a Gaussian random vector and $\{\Psi_j \}$ are the Hermite polynomials. However, alternative representations using different types of random variables and polynomials are also available using generalized PC expansion presented in~\cite{xiu2003modeling}.

The multidimensional polynomial chaoses up to the second order are given as \cite{ghanemSFEM1991} 
\begin{align}
\Psi_0 ( \xi ) & = 1 \ , \\ \nonumber
\Psi_1 ( \xi_{i_1} ) & = \xi_{i_1} \ , \\ \nonumber
\Psi_2 ( \xi_{i_1}, \xi_{i_2} ) & = \xi_{i_1} \xi_{i_2} - \delta_{{i_1} {i_2}}\ ,
\end{align}
where $\delta_{{i_1} {i_2}}$ denotes the Kronecker delta function defined as
\begin{equation}
\delta_{{i_1} {i_2}} = \left\{ \begin{array}{lll}
1 & {\rm if} & {i_1}={i_2} \\
0 & {\rm if} & {i_1} \neq {i_2}
\end{array} \right. \ .
\end{equation}
Note that, $ \{ {\Psi}_i \} $ are orthogonal in the statistical sense, i.e., their inner product $ \left< \Psi_{i_1} \Psi_{i_2} \right>$ is zero for $i_1 \neq i_2$.
For example, the second order PCE of $\alpha (\textbf{\textit{x}},  {\bxi})$ with three random variables $\left\{ \xi_1,\xi_2, \xi_3 \right\}$ is expanded as~\cite{ghanemSFEM1991},
\begin{align}\label{eq:2Dhpoly}
\alpha (\textbf{\textit{x}},  {\bxi}) &= \alpha_0 \Psi_0 + \alpha_1 \Psi_1  + \alpha_2 \Psi_2 + \alpha_3 \Psi_3 + \alpha_4 \Psi_4 + \alpha_5 \Psi_5
 \\ \nonumber
& \ \ \ \ \  + \alpha_6 \Psi_6 + \alpha_7 \Psi_7 + \alpha_8 \Psi_8 + \alpha_9 \Psi_9,  \\
 &= \alpha_0 + \alpha_1 \xi_1 + \alpha_2 \xi_2 + \alpha_3 \xi_3
 + \alpha_4 ( \xi_1^2\ -\ 1)  + \alpha_5 \xi_1 \xi_2 \\ \nonumber & \ \ \ \ \
 + \alpha_6 ( \xi_2^2\ -\ 1) + \alpha_7 \xi_1 \xi_3 + \alpha_8 (\xi_3^2\ -\ 1)
+ \alpha_9 \xi_2 \xi_3.
\end{align}

The explicit expressions for the polynomials used in \Cref{eq:2Dhpoly}  are shown in \Cref{table:hpoly}. The number of terms $P_u$ required in a PCE, with order $p$ and dimension $L$ can be obtained as~\cite{ghanemSFEM1991},
\begin{equation}\label{apn:nPCE}
 P_u = \frac{( L+p ) !}{L! p!}  - 1\ .
\end{equation}

\begin{table}[htbp]
\onehalfspacing
\begin{center}
\begin{threeparttable}
\caption{Polynomial chaoses and their variances for a second-order and three-dimensional PCE~\cite{ghanemSFEM1991}}
\begin{tabular}{| c | c | c | c |} \hline
$j^{th}$ PC term & Order of the expansion & $\Psi_{j}$&${\langle} \Psi_{j}^2 {\rangle}$ \\ \hline
0 & 0 & 1 & 1 \\ \hline
1 & 1 & $\xi_1$ & 1 \\
2 &  & $\xi_1$ & 1 \\
3 &  & $\xi_3$ & 1 \\ \hline
4 & 2 & $\xi_1^2\ -\ 1$ & 2 \\
5 &  & $\xi_1 \xi_2$ & 1 \\
6 &  & $\xi_2^2\ -\ 1$ & 2 \\
7 &  & $\xi_1 \xi_3$ & 1 \\
8 &  & $\xi_3^2\ -\ 1$ & 2 \\
9 &  & $\xi_2 \xi_3$ & 1\\ \hline
\end{tabular}
\label{table:hpoly}
\end{threeparttable}
\end{center}
\end{table}

For the numerical implementation of PCE, one can generalized the evaluation of multidimensional polynomials $ \{ {\Psi}_i \} $ using one-dimensional polynomials $ \{ {\psi}_i \} $ and multi-index ${m^j_i}$ defined in~\cite{le2010spectral}.
This approach is quite useful in the automation of PCE basis function evaluation and calculation of their moments for high-dimensional PC expansions.
For demonstration, consider evaluations of the $L$-dimensional polynomial chaoses and their moments using one-dimensional polynomials.
The form of the one-dimensional Hermite polynomials are given below.
\begin{align}
{\psi}_0  &=   1, \nonumber \\
{\psi}_1  &=   \xi , \nonumber \\
{\psi}_2  &=   {\xi }^2-1, \nonumber \\
{\psi}_n  &=   \xi{\psi}_{n-1}-(n-1){\psi}_{n-2.}.
\label{eq:1dhermite}
\end{align}
The $L$-dimensional polynomial chaoses can be obtained from~\cite{le2010spectral}:
\begin{align}\label{eq:mhermite}
  {\Psi}_{j}({\xi_1,\xi_2,\dots,\xi_L}) = \prod_{i=1}^L \psi_{m^j_i}(\xi_i)
\end{align}
where ${m^j_i}$ denotes multi-index. The code adapted from UQTk~\cite{debusschere2013uqtk} is employed in this thesis to get the multi-index (refer to~\cite{le2010spectral} for further details on multi-index definition and construction). A snippet of the Matlab code used for the evaluation of multidimensional Hermite polynomials in \Cref{eq:mhermite} is given in Listing~\ref{list:multiPoly}.
Similar procedure can be employed to evaluate the moments of multidimensional polynomials $ \{ {\Psi}_i \} $ using moments of one-dimensional polynomials and multi-index~\cite{le2010spectral}.
The moment of multidimensional polynomials $ \{ {\Psi}_i \} $ of order $p$ and dimension $L$ can be obtained using one-dimensional polynomials $ \{ {\psi}_i \} $ and multi-index $m^j_i$ as:
\begin{align}\label{eq:moment}
\left< \prod_{n=1}^p \Psi_{j_{n}} \right>_L = \prod_{i=1}^L \left< \prod_{n=1}^p    \psi_{m^{j_n}_i}\right>_1
\end{align}
The code adapted from UQTk~\cite{debusschere2013uqtk} is employed to evaluate moments of multidimensional polynomials. The Matlab code snippet to evaluate the moments of $L$-dimensional Hermite polynomials using moments of one-dimensional polynomials is given in Listing~\ref{list:cijk}.
Note that the direct evaluation of moments of multidimensional polynomials by solving multidimensional integral is computationally expensive, especially for the high-dimensional cases.

\subsection{Spectral Representation of Lognormal Stochastic Process using PCE}\label{apn:lnp}
The PCE of a lognormal stochastic process $l(\textbf{\textit{x}},\theta)$, obtained by exponential of a Gaussian process $g(\textbf{\textit{x}},\theta)$, with a covariance function $C_{\alpha\alpha}$ and variance $\sigma^2$ defined over a given domain, for instance, as shown in \Cref{eq:apn_exp2d},  
\begin{equation}\label{eq:lxgx}
  l(\textbf{\textit{x}},\theta) = \mathrm{exp} \left[ g(\textbf{\textit{x}},\theta) \right].
\end{equation}
The underlying Gaussian process $g(\textbf{\textit{x}},\theta)$ is characterized by using a truncated KLE with $L$ random variables as follows, 
\begin{equation}\label{eq:gKLE}
g(\textbf{\textit{x}},  {\theta}) =  g_0(\textbf{\textit{x}}) + \sum^L_{j=1} g_j(\textbf{\textit{x}})  \ {\xi_j(\theta)},
\end{equation}
where $g_0(\textbf{\textit{x}})$ is the mean and $g_j(\textbf{\textit{x}})  =  \sqrt{{\lambda}_j} \ f_j(\textbf{\textit{x}})$ with $ \{ \lambda_j \} $ and $ \{ f_j \} $ denoting eigenvalues and eigenvectors respectively as defined in~\Cref{apn:kle}.
The lognormal stochastic process in \Cref{eq:lxgx} can be rewritten using \Cref{eq:gKLE} as
\begin{equation}\label{eq:lxgx2}
  l(\textbf{\textit{x}},\theta) = \mathrm{exp} \left[  g_0(\textbf{\textit{x}}) + \sum^L_{j=1} g_j(\textbf{\textit{x}})  \ {\xi_j(\theta)} \right].
\end{equation}
The lognormal process $l(\textbf{\textit{x}},\theta)$ can be expanded using PCE as follows,
\begin{equation}\label{eq:lpce}
l(\textbf{\textit{x}},  {\theta}) =   \sum^{P_\alpha}_{i=0} l_i(\textbf{\textit{x}})  \ {\Psi_i(\bxi)},
\end{equation}
where $P_\alpha$ is number of PCE terms obtained by using \Cref{apn:nPCE} and $\{ l_i(\textbf{\textit{x}}) \}_{i=0}^{P_\alpha}$ are the PCE coefficients of the lognormal process $l(\textbf{\textit{x}},\theta)$.

Performing Galerkin projection, $l_i(\textbf{\textit{x}})$ can be obtained as~\cite{ghanemSFEM1991,ghanem1999ingredients}
\begin{equation}\label{eq:lpce_coe}
l_i(\textbf{\textit{x}}) = \frac{\left< l(\textbf{\textit{x}},  {\theta}) {\Psi_i(\bxi)} \right>}{\left< {\Psi^2_i(\bxi)} \right>}.
\end{equation}
The denominator in \Cref{eq:lpce_coe} can be evaluated analytically beforehand, for instance see \Cref{table:hpoly}.
The numerator in \Cref{eq:lpce_coe} can be expressed as an  integral~\cite{ghanem1999ingredients}
\begin{equation}\label{eq:lpce_int}
\left< l(\textbf{\textit{x}},  {\theta}) {\Psi_i(\bxi)} \right>
= \int_{-\infty}^{+\infty} \mathrm{exp} \left[  g_0(\textbf{\textit{x}}) + \sum^L_{j=1} g_j(\textbf{\textit{x}}) \xi_j \right] {\Psi_i(\bxi)}   \ \mathrm{exp} \left[ -\frac{1}{2} \sum^L_{j=1} \xi_j^2 \right] d\bxi.
\end{equation}
\Cref{eq:lpce_int} can be simplified to~\cite{subber2012PhDTh}
\begin{equation}\label{eq:lpce_num}
\left< l(\textbf{\textit{x}},  {\theta})  {\Psi_i(\bxi)} \right>
= \mathrm{exp} \left[  g_0(\textbf{\textit{x}}) + \frac{1}{2}\sum^L_{j=1} g_j^2(\textbf{\textit{x}}) \right]  \left< {\Psi_j(\boldsymbol{\eta})} \right> .
\end{equation}
where $\boldsymbol{\eta}_j = \bxi_j - g_j(\textbf{\textit{x}})$.
\Cref{eq:lpce_num} can be rewritten in more concise form as
\begin{equation}\label{eq:lpce_num2}
\left< l(\textbf{\textit{x}},  {\theta}) , {\Psi_i(\bxi)} \right>
= l_0(\textbf{\textit{x}}) \left< {\Psi_j(\boldsymbol{\eta})} \right>.
\end{equation}
where $l_0(\textbf{\textit{x}})$ represents the mean of the lognormal process $l(\textbf{\textit{x}},  {\theta})$
\begin{equation}
l_0(\textbf{\textit{x}}) = \mathrm{exp} \left[  g_0(\textbf{\textit{x}}) + \frac{1}{2}\sum^L_{j=1} g_j^2(\textbf{\textit{x}}) \right].
\end{equation}
As the number of KLE terms, $L$ tends to $\infty$, the mean of $l(\textbf{\textit{x}},  {\theta})$ converges to~\cite{ghanem1999ingredients}
\begin{equation}
l_0(\textbf{\textit{x}}) = \mathrm{exp} \left[  g_0(\textbf{\textit{x}}) + \frac{1}{2}\sigma^2 \right].
\end{equation}

Using \Cref{eq:lpce_coe} and \Cref{eq:lpce_num}, the PCE of the lognormal stochastic process $l(\textbf{\textit{x}},{\theta})$ can be written as
\begin{equation}\label{eq:apn_logPCE}
l(\textbf{\textit{x}},\theta) = l_0(\textbf{\textit{x}}) \sum^{P_{\alpha}}_{i=0} \frac{\left< {\Psi_i(\boldsymbol{\eta})} \right>}{\left<  \Psi^2_i(\bxi) \right>}  \Psi_i(\bxi).
\end{equation}

The $\left< {\Psi_j(\boldsymbol{\eta})} \right>$ represents the expectation of the PC basis around the coefficients $g_i(\textbf{\textit{x}})$.
These expectations can be evaluated analytically. For instance, see \Cref{table:logNoraml} showing variance $\left<  \Psi^2_i(\bxi) \right>$ and expectations $\left< {\Psi_i(\boldsymbol{\eta})} \right>$ for the second-order and three-dimensional PC basis functions.
\begin{table}[htbp]
\onehalfspacing
\begin{center}
\begin{threeparttable}
\caption{The expectation and variance of second-order and three-dimensional PC basis~\cite{ghanem1999ingredients}.}
\begin{tabular}{| c | c | c | c |} \hline
$\Psi_i(\bxi)$ & $ \left<  \Psi^2_i(\bxi) \right> $ & $\left< {\Psi_j(\boldsymbol{\eta})} \right>$ \\ \hline
$\xi_1$ & 1 & $g_1(\textbf{\textit{x}})$\\
$\xi_1$ & 1 & $g_2(\textbf{\textit{x}})$ \\
$\xi_3$ & 1 & $g_i(\textbf{\textit{x}})$ \\ \hline
$\xi_1^2\ -\ 1$ & 2 & $g^2_1(\textbf{\textit{x}})$ \\
$\xi_1 \xi_2$ & 1 & $g_1(\textbf{\textit{x}}) g_2(\textbf{\textit{x}})$ \\
$\xi_2^2\ -\ 1$ & 2 & $g^2_2 (\textbf{\textit{x}})$ \\
$\xi_1 \xi_3$ & 1 & $g_1(\textbf{\textit{x}}) g_3(\textbf{\textit{x}})$ \\
$\xi_3^2\ -\ 1$ & 2 & $g^2_3(\textbf{\textit{x}})$ \\
$\xi_2 \xi_3$ & 1 & $g_2(\textbf{\textit{x}}) g_3(\textbf{\textit{x}})$ \\ \hline
\end{tabular}
\label{table:logNoraml}
\end{threeparttable}
\end{center}
\end{table}

For further simplification, using $L=3$ and the second order ($p=2$)  expansion leading to the number of PCE terms $P_\alpha = 9$, \Cref{eq:apn_logPCE} can be expand as (note $\textbf{\textit{x}}$ of $g_i$ is dropped for notational convenience),
\begin{align}\label{eq:logPCE3}
l(\textbf{\textit{x}},\theta) &= l_0 \Big( 1 + \xi_1 g_1 + \xi_2 g_2 + \xi_3 g_3 + (\xi_1^2-1) \frac{g_1^2}{2} \nonumber \\
&+  (\xi_1 \xi_2) g_1 g_2 + (\xi_2^2-1) \frac{g_2^2}{2} +( \xi_1 \xi_3) g_1 g_3 + (\xi_3^2-1) \frac{g_3^2}{2} + (\xi_2 \xi_3) g_2 g_3 \Big)
\end{align}

For illustration, consider a simplest case where $g(\textbf{\textit{x}},\theta)$ is characterized by using a Gaussian random variable $\xi$ with the mean $\mu_g$ and variance $\sigma^2_g$.
The lognormal random variable $l(\textbf{\textit{x}},\theta)$ can be obtained using procedure outlined above~\cite{ghanem1999ingredients},
\begin{align}\label{eq:logPCE2}
l(\textbf{\textit{x}},\theta) &= \mu_l \sum^L_{i=0} \frac{\sigma^i_g}{i!}  \Psi_i(\xi), \nonumber \\
 &= \mu_l \left( 1 + \frac{\sigma_g}{1!} \xi + \frac{\sigma^2_g}{2!} (\xi^2-1)  + \dots \right),
\end{align}
where $\mu_l =\mathrm{exp} \left[  \mu_g + \frac{1}{2}\sigma^2_g \right]$ is the mean of the lognormal random variable. 

\section{Spectral Stochastic Finite Element Method}\label{sec:ssfem}

Consider a two-dimensional steady-state flow through random media with a spatially varying non-Gaussian diffusion coefficient $c_d$. The flow is modeled by a two-dimensional stochastic diffusion equation. This leads to a Poisson problem defined by a linear elliptic stochastic PDE 
as defined below:
\begin{align}\label{eq:sspe}
-\nabla \ \cdotp \big( \ c_d(\textbf{\textit{x}},\theta) \ \ \nabla u(\textbf{\textit{x}}, \theta) \ \big) &= {F}(\textbf{\textit{x}}),     \ \ \  \ \ \ \   \mathcal{D}\times \Omega, \\
u(\textbf{\textit{x}}, \theta) &= 0,    \ \ \  \ \ \ \  \ \ \ \ \  \partial \mathcal{D}\times \Omega,
\end{align}
where $\nabla$ denotes the gradient which represents the differential operator with respect to the spatial variables $\textbf{\textit{x}}$, $u$ is the solution process,  $\theta$ is an element in the sample space $\Omega$ defined by the probability space ($\Omega, \mathcal{F}, \mathcal{P})$.
For the sake of convenience, $F(\textbf{\textit{x}})$ is modeled as a deterministic source term. However, the methodology presented herein can be easily extended to stochastic source function $F(\textbf{\textit{x}},\theta)$.


The finite element discretization with $N$ nodes in the spatial domain leads to a system of linear equations with random coefficients $\theta$ denotes stochasticity~\cite{ghanemSFEM1991}
\begin{equation}\label{eq:dFEM}
{\mathrm{\bf{A}(\theta) \bf{u}(\theta) = \bf{f}}},
\end{equation}
where ${\mathrm{\bf{A}}(\theta)}$ is the random or stochastic system matrix, $ \bf{u}(\theta)$ is the stochastic response vector and $\bf{f}$ is the deterministic source vector. 

The above system can be solved for the mean or any sample of the stochastic parameter $c_d(\textbf{\textit{x}},\theta)$ by using any deterministic FEM solver. For demonstration we employed FEniCS general purpose deterministic FEM packages~\cite{fenics/dolfin17}. 

\begin{figure}[htbp]
 \centering
 \subfloat[Two-dimensional mesh \label{subfig:mesh}]{%
  \includegraphics[width=0.42\textwidth,height=0.27\textheight]{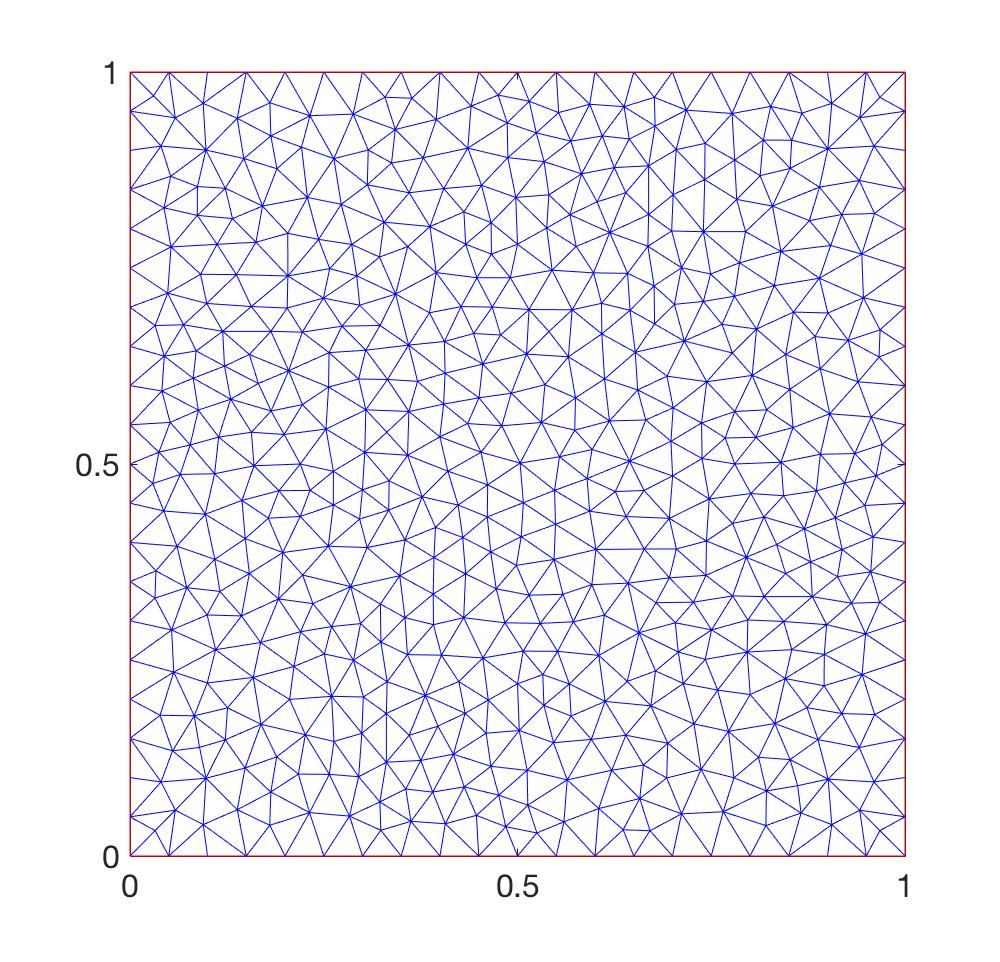}}
  \subfloat[The mean solution field \label{subfig:mesh}]{%
  \includegraphics[width=0.5\textwidth,height=0.27\textheight]{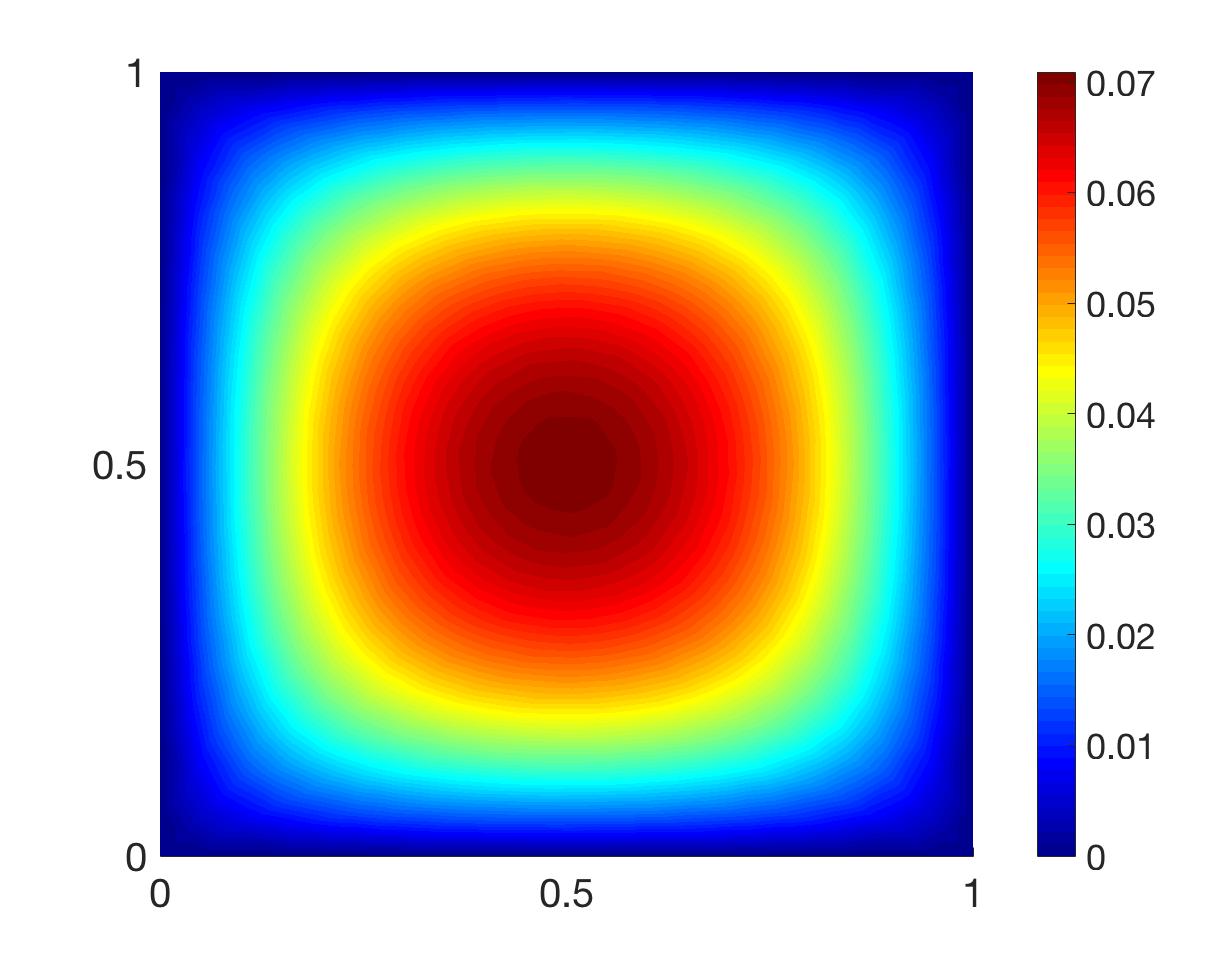}}
  \caption{Finite element mesh and the solution field at the mean value of stochastic system parameters.}
  \label{fig:mean_SSFEM}
\end{figure}

Consider a unit square domain discretized using unstructured finite element mesh with $600$ nodes and $1200$ elements as shown in \Cref{fig:mean_SSFEM} (left). The numerical simulations are performed using the unit source term. The solution field for the mean value of stochastic parameter is shown in~\Cref{fig:mean_SSFEM}.
The corresponding FEniCS based python code snippet is shown in Listing~\ref{list:det_fem}.

To solve the above PDE using SSFEM the stochastic system matrix $\mathrm{\bf{A}(\theta)}$ and the stochastic solution process $\bf{u}(\theta)$ in \Cref{eq:dFEM} are approximated using the polynomial chaos expansions as~\cite{ghanemSFEM1991},
\begin{align}
{\mathrm{\bf{A}(\theta)}} \approx  \sum^{P_{\textsc{a}}}_{i=0} {\hat{{\mathrm{\bf{A}}}}_i} \ {\Psi_i(\bxi)},  \ \ \ \ \  {\mathrm{\bf{u}(\theta)}} \approx  \sum^{P_{u}}_{j=0} {\mathrm{\hat{\bf{u}}}_j}  \ {\Psi_j(\bxi)},
\label{eq:apn_pceinout}
\end{align}
where $\hat{{\mathrm{\bf{A}}}}_i$'s are the PCE coefficients of the random system matrix, $\hat{{\bf{u}}}_j$'s are the PCE coefficients of the solution process and $\Psi_j$'s are the multidimensional polynomials obtained as a function of $L$ random variables ${\bxi} = \{ \xi_{1}, \xi_{2}, \dots, \xi_{L} \}$.
$P_{\textsc{a}}$ and $P_{u}$ are the numbers of PCE terms required to express the stochastic system matrix and the solution process, respectively~\cite{ghanemSFEM1991}.
The PCE coefficients ($\hat{{\mathrm{\bf{A}}}}_i$'s) 
are computed using a lognormal diffusion coefficient obtained from a underlying Gaussian process expanded using KLE~\cite{ghanemSFEM1991,pellissetti2000iterative} (refer to \Cref{apn:kle}). Therefore, in the SSFEM approaches, the primary goal is to estimate the PCE coefficients $({\mathrm{\hat{\bf{u}}}_j})$ of the solution process ${\mathrm{\bf{u}(\theta)}}$. In the subsequent sections, we illustrate how to solve this PDE using intrusive and non-intrusive SSFEM. 

\subsection{Intrusive SSFEM}\label{sec:intrusive}
In the intrusive SSFEM, the PCE of the system matrix with random coefficients ${\mathrm{\bf{A}(\theta)}}$ and the solution process ${\mathrm{\bf{u}(\theta)}}$, presented in \Cref{eq:apn_pceinout} are directly substituted into the finite element discretization of stochastic PDE given in \Cref{eq:dFEM} leading to~\cite{ghanemSFEM1991}
\begin{equation}
\epsilon = \sum^{P_{\textsc{a}}}_{i=0} {\hat{{\mathrm{\bf{A}}}}_i} {\Psi_i({\pmb{\xi}})}
\ \sum^{P_{u}}_{j=0} {\mathrm{\hat{\bf{u}}}_j} {\Psi_j({\pmb{\xi}})} - {\mathrm{\bf{f}}} \neq 0,
\label{eq:skuf}
\end{equation}
where $\epsilon$ is the random residual.

Performing Galerkin projection, i.e., multiplying both sides of the above equation by ${\Psi_k({\pmb{\xi}})}$ with $k = 0,..,P_{u}$ and taking expectation both sides results in the following system of coupled equations~\cite{ghanemSFEM1991},
\begin{equation}
  \left< \epsilon , {\Psi_k({\pmb{\xi}})} \right> = 0 \ , \ \ \ \ \   k = 0,1,\dots,P_{u}
\end{equation}
\begin{equation}
\sum^{P_{u}}_{j=0} \sum^{P_{\textsc{a}}}_{i=0} \left< {\Psi_i({\pmb{\xi}})} {\Psi_j({\pmb{\xi}})} {\Psi_k({\pmb{\xi}})} \right>  {\hat{{\mathrm{\bf{A}}}}_i}  {\mathrm{\hat{\bf{u}}}_j}    =  \left< {\mathrm{\bf{f}}} \ {\Psi_k({\pmb{\xi}})} \right> , \ \ k = 0,1,\dots,P_{u}.
\label{eq:gpkuf}
\end{equation}
For notational convenience, we rewrite \Cref{eq:gpkuf} using $\left< {\Psi_i({\pmb{\xi}})} {\Psi_j({\pmb{\xi}})} {\Psi_k({\pmb{\xi}})} \right> = \mathbfcal{C}_{ijk}$  and $ \left< {\mathrm{\bf{f}}} \ {\Psi_k({\pmb{\xi}})} \right> = f_{k}$ leading to
\begin{equation}
\sum^{P_{u}}_{j=0} \sum^{P_{\textsc{a}}}_{i=0}  \mathbfcal{C}_{ijk} {\hat{{\mathrm{\bf{A}}}}_i} \ {\mathrm{\hat{\bf{u}}}_k}   =  f_{k} \ , \ \ \ \ \   k = 0,1,\dots,P_{u}.
\label{eq:cijkkuf}
\end{equation}
For concise representation, the following notation is used,
\begin{equation}\label{eq:sFEMatrix}
A_{jk} =  \sum^{P_{\textsc{a}}}_{i=0}  \mathbfcal{C}_{ijk} {\hat{{\mathrm{\bf{A}}}}_i},
\end{equation}
Thus, \Cref{eq:cijkkuf} can be further simplified as,
\begin{equation}
\sum^{P_{u}}_{j=0}  A_{jk}  \ {\mathrm{\hat{\bf{u}}}_k}    = f_k \ ,  \ \ \ \ \   k = 0,1,\dots,P_{u}.
\label{eq:Gakuf}
\end{equation}
Equivalently \Cref{eq:Gakuf} can be written as,
\begin{equation}\label{eq:sFEM}
\mathcal{[A] \{U\} = \{F\}},
\end{equation}
where $\mathcal{A}$ in \Cref{eq:sFEM} is the system matrix,
$\mathcal{U}$ is the vector of PCE coefficients of the solution process and
$\mathcal{F}$ is the corresponding right hand side vector arising in the setting of intrusive SSFEM. 
The size of the system matrix $\mathcal{A}$ is $(N\times P_{u} ,  N\times P_{u})$, where $N$ is the number of degree-of-freedom related to the finite element mesh resolution and $P_{u}$ is the number of PCE terms used in the representation of the solution process.
Note that, $P_{u}$ is a function of stochastic dimension $L$ and order $p_{u}$ of the PCE.
From \Cref{eq:sFEMatrix} and \eqref{eq:Gakuf} it can be noted that, each block of the system matrix $\mathcal{A}$ is denoted by $A_{jk}$ (a sub matrix of size $(N \times N))$ and it can be computed from the set of deterministic finite element matrices ${\hat{{\mathrm{\bf{A}}}}_i}$. 

In this article, the system matrix assembly procedure, i.e., implementation of~\Cref{eq:sFEMatrix} is performed by
employing deterministic finite element assembly routines imported from the FEniCS general purpose FEM package.
The procedure for stochastic system matrix assembled is outline in \Cref{alg:smap_local}. This procedure employs deterministic, element-level, FEniCS-based assemble routines. It is designed to reduce the number of call to the deterministic assembly routine.

\begin{figure}[htbp!]
 \centering
  \subfloat[Mean]{%
  \includegraphics[width=0.46\textwidth,height=0.27\textheight]{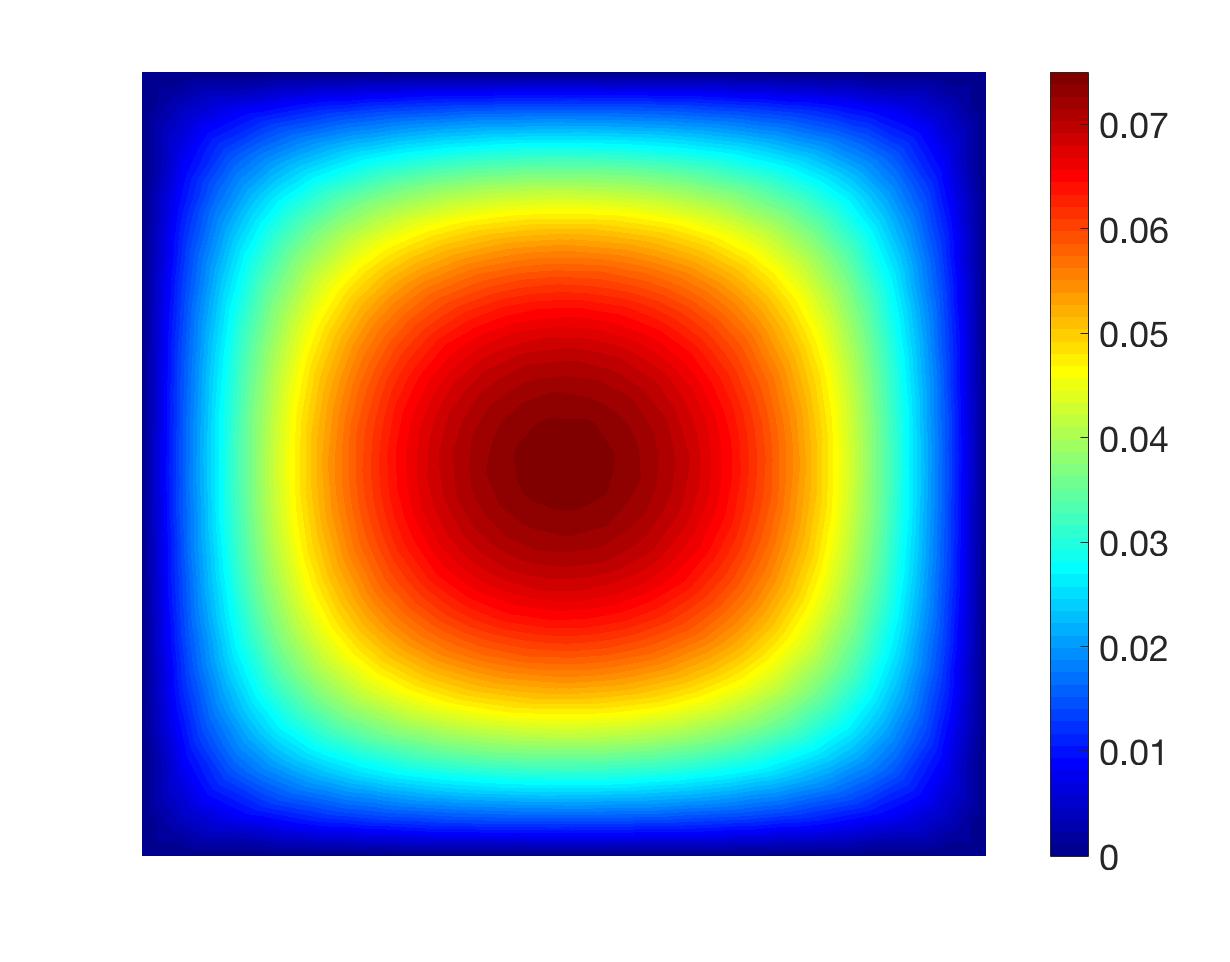}}
  \subfloat[Standard deviation]{%
  \includegraphics[width=0.46\textwidth,height=0.27\textheight]{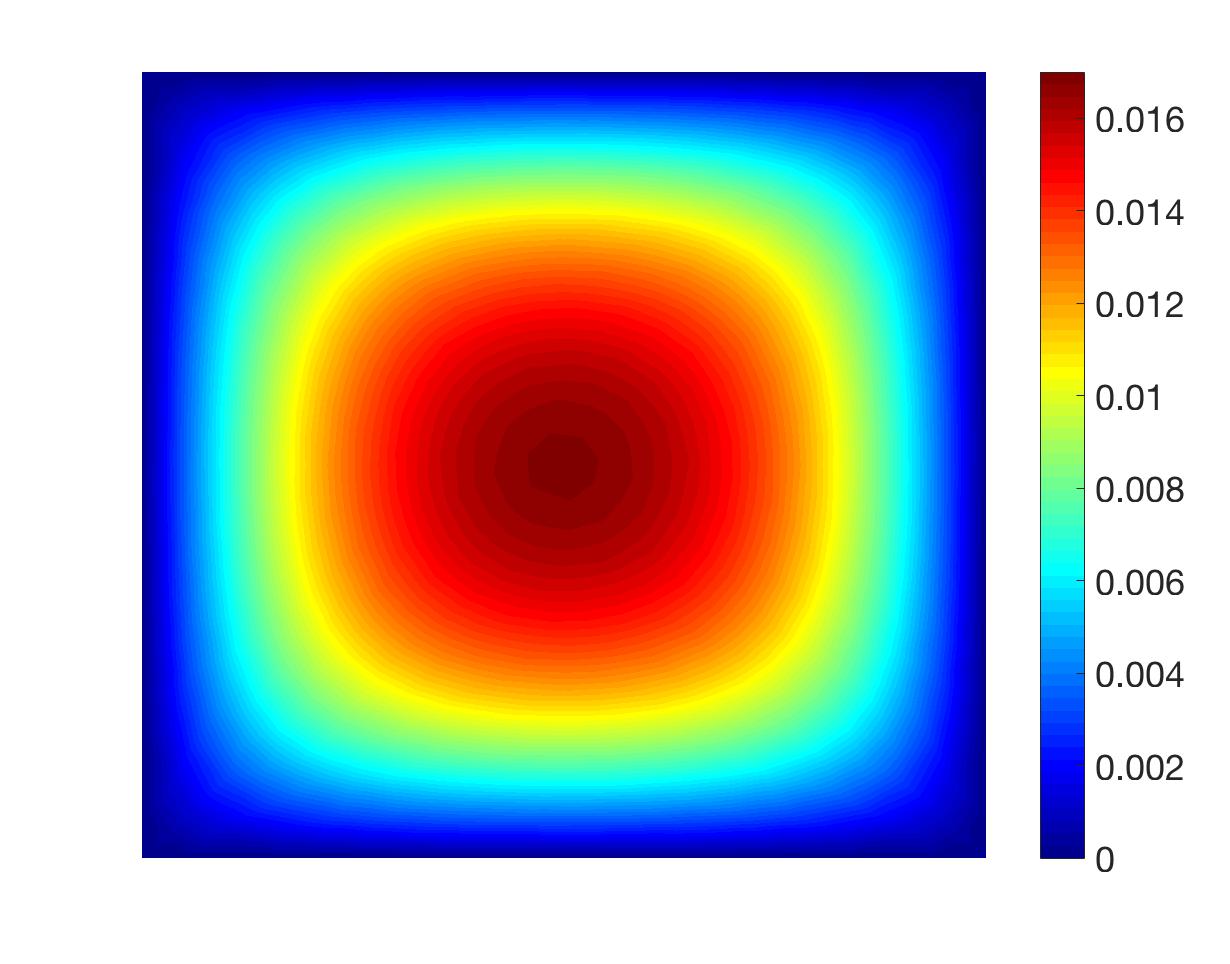}}
  
   \subfloat[${\mathrm{\hat{\bf{u}}}_4}$ \label{subfig:in4}]{
  \includegraphics[width=0.46\textwidth,height=0.27\textheight]{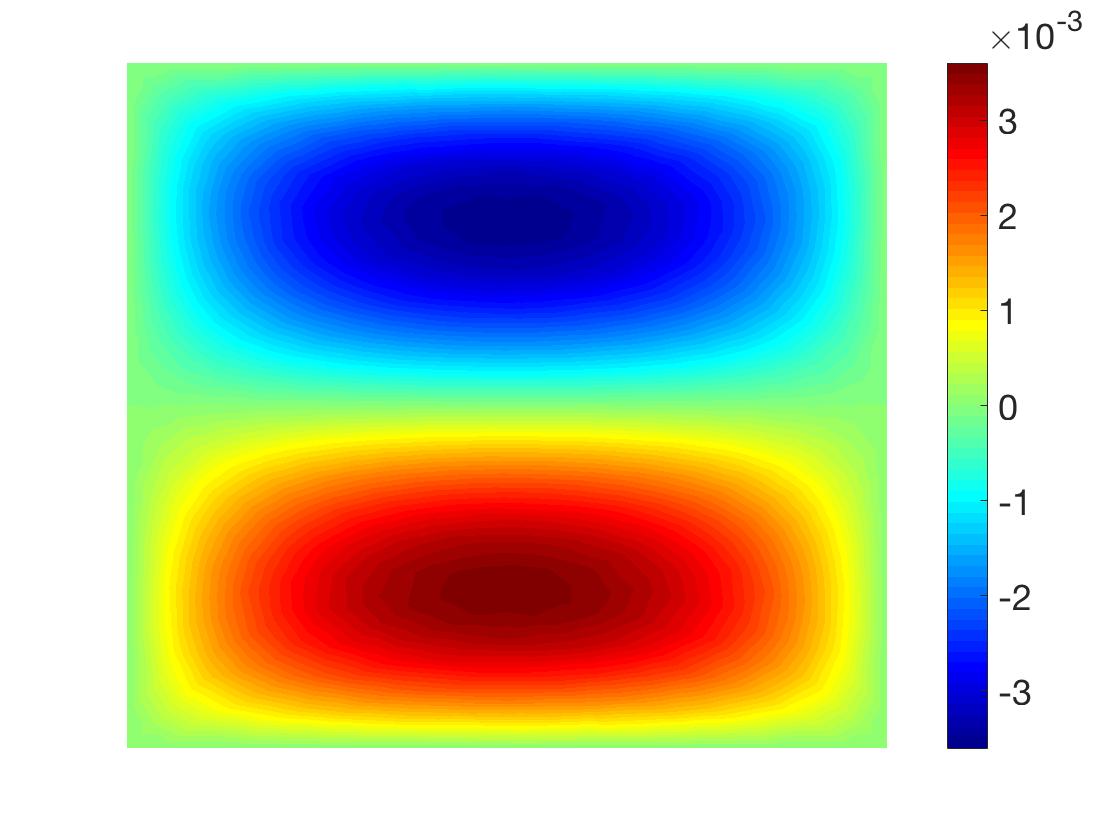}}
  \subfloat[${\mathrm{\hat{\bf{u}}}_5}$ \label{subfig:nisp}]{
  \includegraphics[width=0.46\textwidth,height=0.27\textheight]{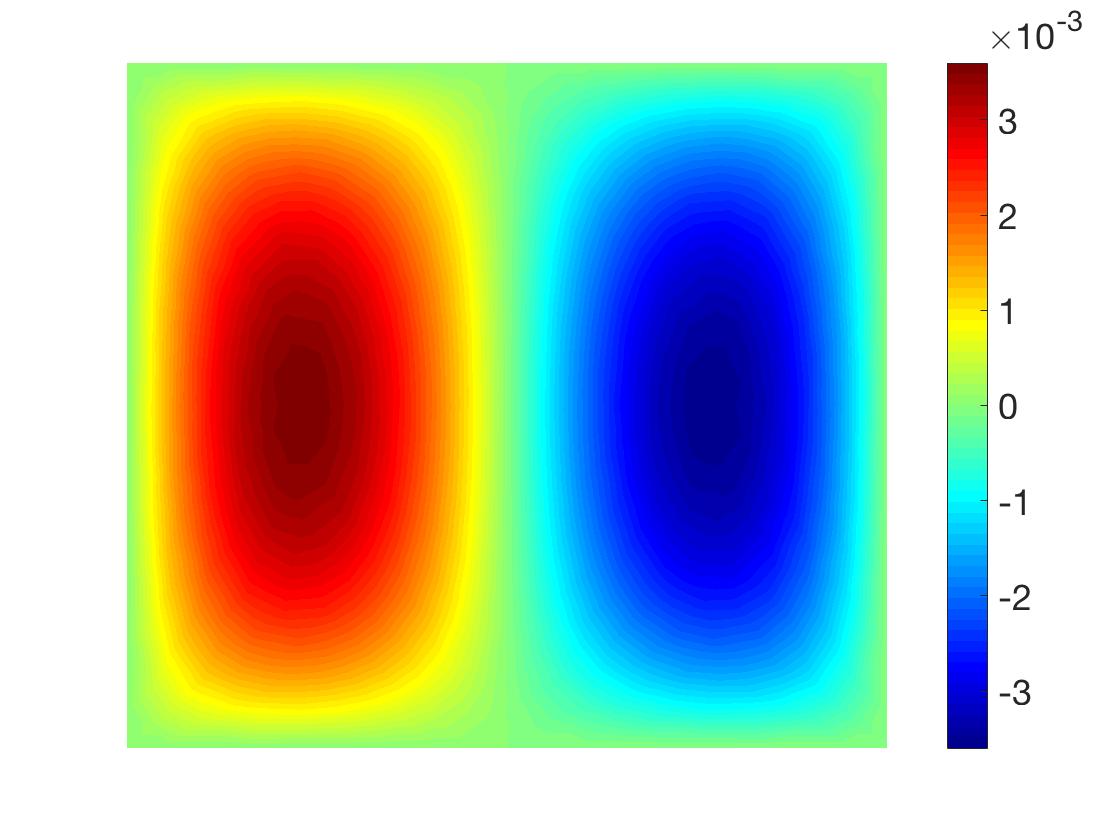}}

  \subfloat[${\mathrm{\hat{\bf{u}}}_9}$ \label{subfig:nisp9}]{
  \includegraphics[width=0.46\textwidth,height=0.27\textheight]{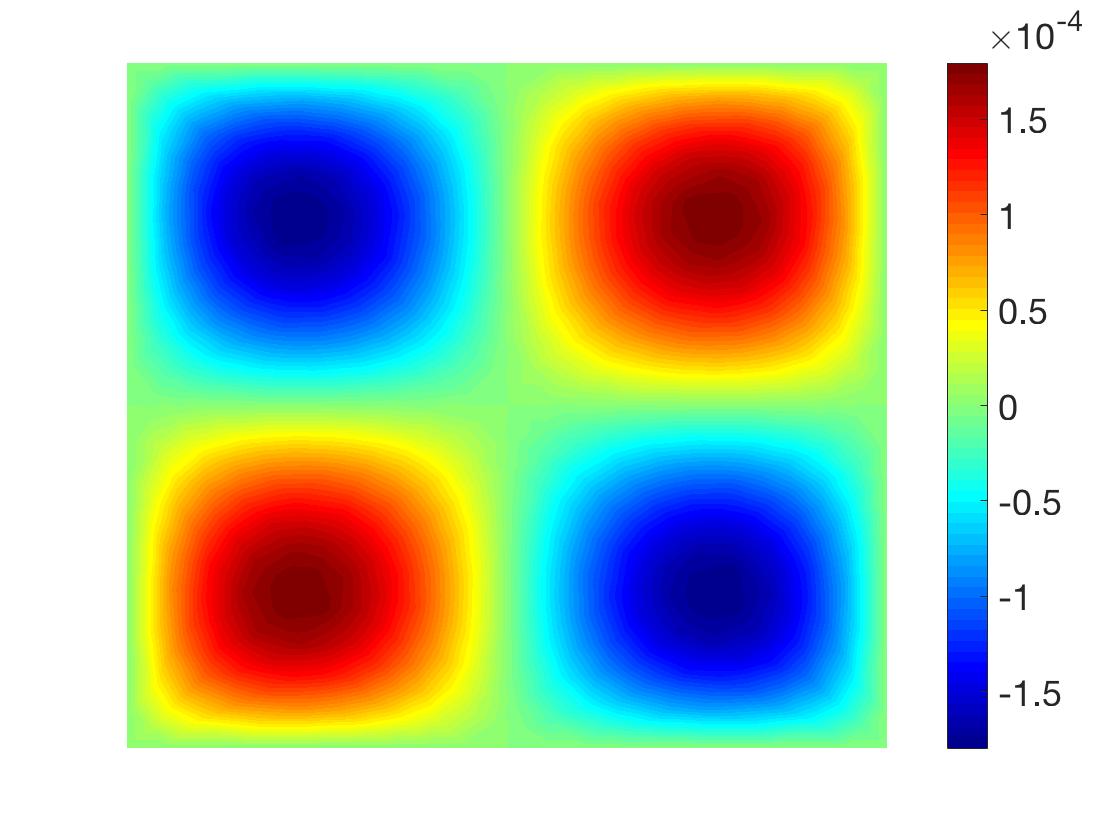}}
  \subfloat[${\mathrm{\hat{\bf{u}}}_{10}}$ \label{subfig:nisp10}]{
  \includegraphics[width=0.46\textwidth,height=0.27\textheight]{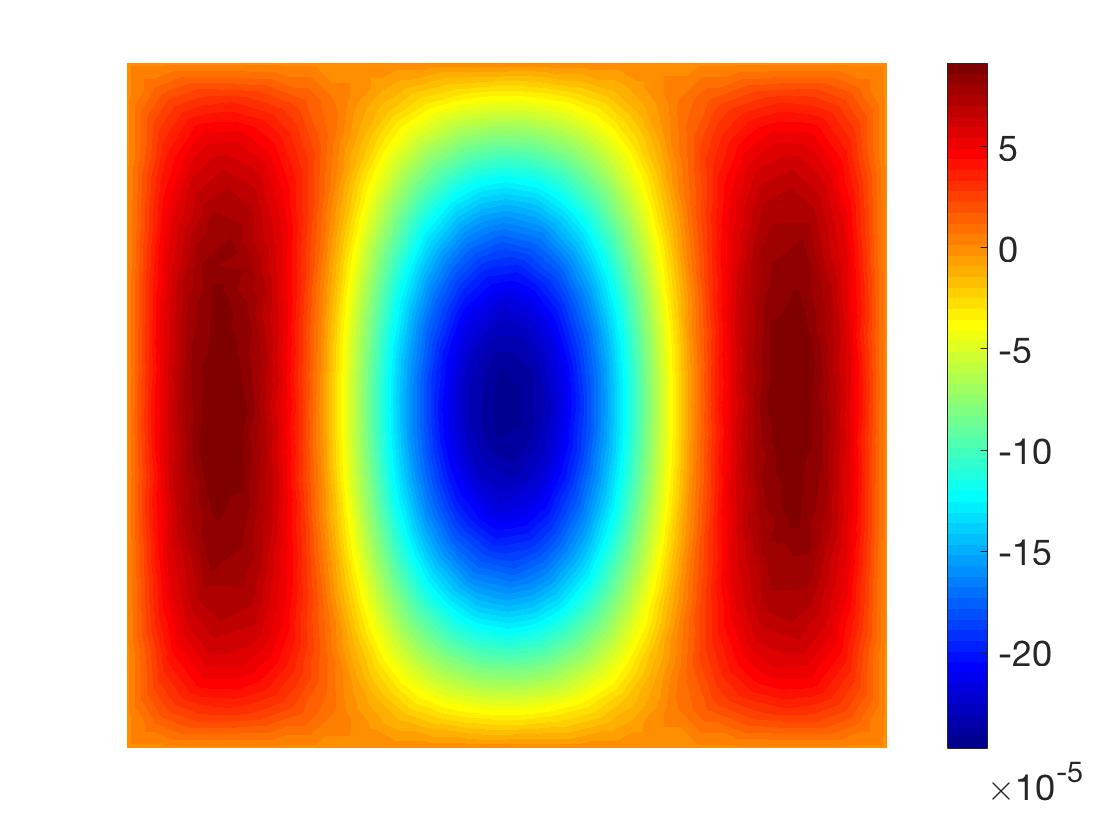}}
  
  \caption{The mean, standard deviation, and a few selected PCE coefficients of the solution field using intrusive SSFEM.}
  \label{fig:intrusiv_sol}
\end{figure}

\alglanguage{pseudocode}
\begin{algorithm}
\caption{: \textcolor{blue}{Stochastic System Matrix Assembly Procedure}}
\label{alg:smap_local}
\begin{algorithmic}[1]
\State {\bf{Input} :}  non-zero : $ijk$ and $C_{ijk}$
\For{$ {i = 0, 1, 2, ..., P_{\textsc{a}}}   $}
\State Call Modified \textcolor{red}{Element Level FEniCS Deterministic Assembly} ( ${\hat{{\mathrm{\bf{[A]}}}}_i}$ )
\For{$ {j = 0, 1, 2, ..., P_{u}}  $}
\For{$ {k = 0, 1, 2, ..., P_{u}}  $}
\If{$ i, j, k == non-zero(ijk) $}
\State \textcolor{red}{Stochastic Matrix Assembly} ( $\mathcal{[A]}_{j,k} \ += C_{ijk}*{\hat{{\mathrm{\bf{[A]}}}}_i}$ )
\EndIf
\EndFor
\EndFor
\EndFor
\State {\bf{Output} :}  $\mathcal{[A]}$
\end{algorithmic}
\end{algorithm}

The code snippet to perform intrusive SSFEM system matrix assembly using FEniCS-based assembly routine is presented in Listing~\ref{list:sto_femAssembly}.
The procedure to call to ${{assemble\_local}}$, an element-level FEniCS assembly routines~\cite{fenics/dolfin17} is outlined in Listing~\ref{list:det_femAssembly}. For each of the input PCE term, the procedure outlined in the Listing~\ref{list:det_femAssembly} is invoked once. The FEniCS-based procedure to define stochastic variational form for steady state diffusion equation defined in~\Cref{eq:sFEM} is outlined in Listing~\ref{list:sto_variation}. The diffusion coefficient $c_d$ is characterized as a lognormal stochastic process. The diffusion coefficient $c_d$ is defined as a FEniCS-based, expression-class~\cite{fenics/dolfin17}, as outlined in Listing~\ref{list:sto_cd}.

The parameters such as $\omega$, $\lambda$, sort-index and multi-index required to perform the SSFEM-system matrix assembly procedure, are need to be computed beforehand. The procedure outlined earlier in \Cref{apn:kle} and \Cref{apn:pce} can be used to calculate $\omega_i$, $\lambda_i$, and sort-index. The multi-index calculation is performed using functions from UQTk.
The procedure to calculate moments of the multidimensional polynomials $C_{ijk}$ and non-zero $i,j$ and $k$ indices, required for the Listing~\ref{list:sto_femAssembly} is outlined in the Listing~\ref{list:cijk}. This procedure is developed by adapting functions from UQTk.

The numerical simulations are performed using FEniCS package for the stochastic PDE defined by \Cref{eq:sspe}, using the following numerical parameters $b=1$, $\sigma=0.3$, $L=3$, $p_u=3$ and $f=1$.
A unit square domain, discretized using $600$ nodes and $1200$ elements 
is used. The mean, standard deviation, and a few selected PCE coefficients of the solution field are plotted in \Cref{fig:intrusiv_sol}.



\subsection{Non-Intrusive SSFEM}\label{sec:nisp}
Performing Galerkin projection onto the PC expansion of solution process given in \Cref{eq:apn_pceinout} and then exploiting orthogonality properties of the basis functions, the PCE coefficients of the solution process can be evaluated as follows~\cite{ghanemSFEM1991,le2010spectral};
\begin{equation}
{\mathrm{\hat{\bf{u}}}_k}   = \frac{\left< {\mathrm{\bf{u}(\theta)}} {\Psi_k({\pmb{\xi}})}  \right>} {\left< {\Psi_j({\pmb{\xi}})} {\Psi_k({\pmb{\xi}})}  \right>} = \frac{1} {\left<{\Psi_k({\pmb{\xi}})^{2}}  \right>} \int_{\Omega} {\mathrm{\bf{u}(\theta)}} {\Psi_k({\pmb{\xi}})} {\mathrm{p}}(\pmb{\xi}) {d{\pmb{\xi}}},
\label{eq:apn_NISPpce}
\end{equation}
where $\left< {\Psi_j({\pmb{\xi}})} {\Psi_k({\pmb{\xi}})}  \right>$ is non-zero only for $j=k$ and it can be obtained analytically beforehand, for instance see \Cref{table:hpoly}. Therefore, the major computational efforts lies in the evaluation of the multidimensional integral in the numerator of \Cref{eq:apn_NISPpce}~\cite{le2010spectral,nobile2008sparse}.

Consider the FE discretization of a stochastic PDE given in \Cref{eq:dFEM}. Using $\{ \pmb{\xi}_{1} ,  \pmb{\xi}_{2} , \dots, \pmb{\xi}_{n_{s}}  \}$ (where $\pmb{\xi}_i$ is a set of $L$ random variables), the following deterministic system is solved at $n_{s}$ sample points using an existing deterministic solver as a black-box.
\begin{equation}\label{eq:apn_nispSamp}
{\mathrm{\bf{A}(\theta)  \bf{u}(\theta) = \bf{f}}}.
\end{equation}

Computed $\bf{u}(\theta)$ at each sample points are used to calculate the PCE coefficients of the solution process using Smolyak sparse grid quadrature~\cite{le2010spectral,debusschere2013uqtk,nobile2008sparse}.

The sparse grid quadrature rule to integrate a multidimensional function $\mathscr{F} = {\mathrm{\bf{u}(\theta)}} {\Psi_k({\pmb{\xi}})}$ in the numerator of \Cref{eq:apn_NISPpce} 
can be constructed using the univariate quadrature rule ${\mathscr{Q}}_{l}^{(1)} \mathscr{F}$ as follows ~\cite{smith2013uncertainty},
\begin{equation}\label{eq:sg1D}
{\mathscr{Q}}_{l}^{(1)} \ {\mathscr{F}} =  \sum_{q=1}^{N_{s}} {\mathscr{F}}(r_{l}^{q})  \ {w_p}^{q}_{l},
\end{equation}
where the subscript $l$ is the level of quadrature and the superscript denotes the dimension $d$ (in this case $d=1$). $N_{s}$ is the number of quadrature points in the sparse grid. 
The sparse grid nodal set for $d=2$ and $l=3$ can be written as~\cite{smith2013uncertainty}
\begin{align}\label{eq:sparseNodes}
  \Theta_{l=3}^{(d=2)} &= \bigcup_{|l'| \leq l+d-1} \Big( \Theta_{l_1}^{(1)}  \times \Theta_{l_2}^{(1)}  \Big) \\ \nonumber
  &= \Big( \Theta_{1}^{(1)}  \times \Theta_{1}^{(1)}  \Big)  \ \ \ \ \ (l_1=1, l_2=1) \\ \nonumber
  &\cup \Big( \Theta_{1}^{(1)}  \times \Theta_{2}^{(1)}  \Big) \cup \Big( \Theta_{2}^{(1)}  \times \Theta_{1}^{(1)}  \Big)
  \\ \nonumber
  &\cup \Big( \Theta_{1}^{(1)}  \times \Theta_{3}^{(1)}  \Big)  \cup \Big( \Theta_{2}^{(1)}  \times \Theta_{2}^{(1)}  \Big)  \cup \Big( \Theta_{3}^{(1)}  \times \Theta_{1}^{(1)}  \Big),
\end{align}
where $l'=(l_1 \dots l_d)$ with $|l'| = \sum_{i=1}^d l_i$
(a specific example to obtain the sparse grid nodal set involving  $\Theta_{1}^{(1)} \times \Theta_{2}^{(1)}$ is given below).
For the implementation of the multidimensional sparse grid the growth rule in one-dimensional quadrature must be defined. In the current implementation, we have used the Gauss-Hermite quadrature rule~\cite{smith2013uncertainty}.
For further simplification, consider the Gauss-Hermite quadrature rule in one dimension with the nodes and weight specified in \Cref{table:sparseGrid}, for different level of quadrature. Using \Cref{eq:sparseNodes} the nodes and weight for $l=3$ and $d=2$ are obtained as shown in \Cref{table:l2d2sparseGrid}.

\begin{table}[htbp]
\centering
\caption{Nodes and weights for Gauss-Hermite quadrature in one dimension and third level.}
\begin{tabular}{|l|c|c|}   
\hline
level $l$    & nodes $\Theta_l^{(1)}$  & weights $\mathcal{W}_l^{(1)}$\\
\hline
$l=1$      & $\{0 \}$   & $\{ 1 \}$       \\
$l=2$      & $\{-1,  1\}$   & $\{0.5, 0.5\}$       \\
$l=3$      & $\{ -1.7321, 0, 1.7321 \}$ & $\{ 0.167, 0.667, 0.167 \}$       \\
\hline
\end{tabular}
\label{table:sparseGrid}
\end{table}

\begin{table}[htbp]
\centering
\caption{Nodes and weights for Gauss-Hermite quadrature in two dimension and third level.}
\begin{tabular}{|l|c|}   
\hline
nodes $\Theta_2^{(2)}$  & weights $\mathcal{W}_2^{(2)}$ \\
\hline
$\{ -1.732, 	0 \}$ & 0.167\\
$\{ -1.0, 	-1.0\}$ & 0.25 \\
$\{ -1.0, 	0 \}$& -0.5 \\
$\{ -1.0,  1.0\}$ & 0.25 \\
$\{ 0,	 -1.732 \}$& 0.167 \\
$\{ 0,	 -1.0\}$ & -0.5 \\
$\{ 0,	 0 \}$& 1.333 \\
$\{ 0, 1.0 \}$& -0.5 \\
$\{ 0,	 1.732\}$ & 0.167 \\
$\{1.0, 	-1.0 \}$& 0.25 \\
$\{1.0, 	0  \}$& -0.5 \\
$\{1.0, 	1.0\}$ & 0.25 \\
$\{1.732, 	0 \}$& 0.167\\
\hline
\end{tabular}
\label{table:l2d2sparseGrid}
\end{table}

The relations in~\Cref{eq:sparseNodes} are used to get the set of nodes in \Cref{table:l2d2sparseGrid}. For instance, the nodes for $\big( \Theta_1^{(1)}  \times \Theta_2^{(1)} \big)$ are obtained by taking tensor produt of $\Theta_1^{(1)} = \{ 0 \}$ and $\Theta_2^{(1)} = \{ -1, 1 \}$, resulting into the following set of nodes = $\big[ \{0, -1\}, \{0, 1\} \big]$ which correspond to the $6^{rd}$ and the $8^{th}$ row in the first column of \Cref{table:l2d2sparseGrid}~\cite{smith2013uncertainty}. 

\begin{figure}[htbp]
 \centering
  \includegraphics[width=0.4\textwidth,height=0.20\textheight]{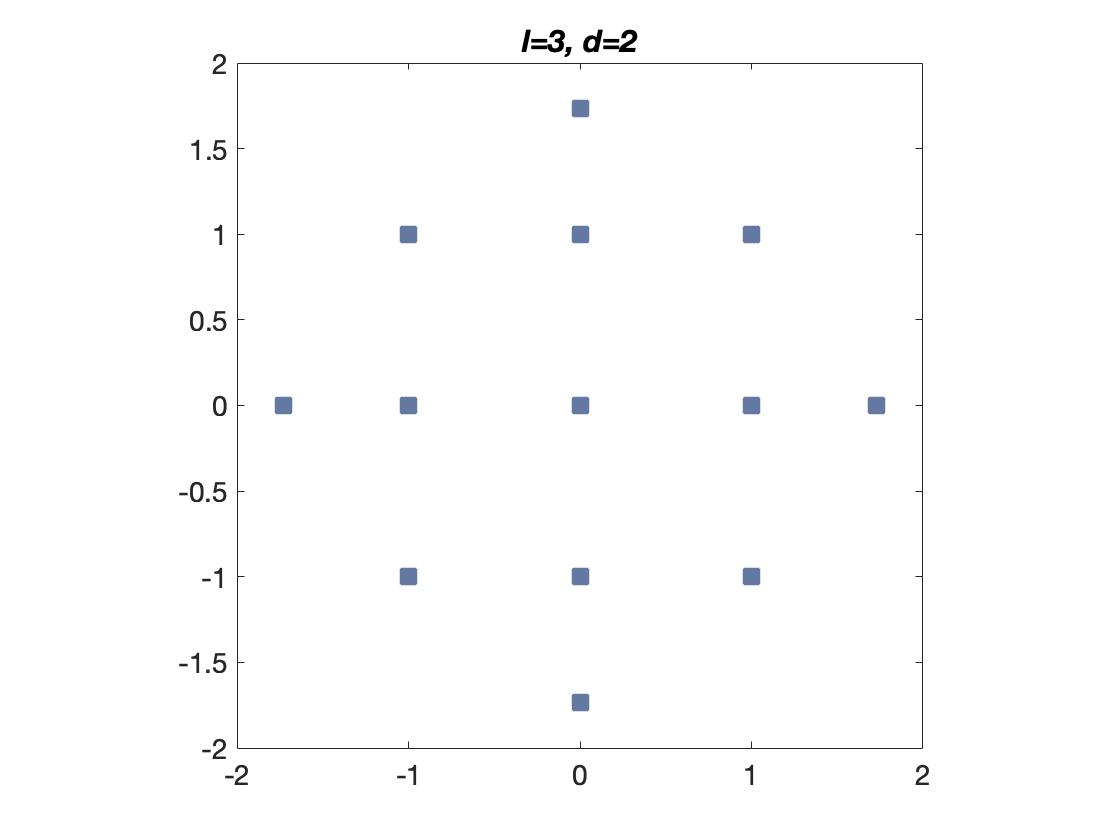}
  \includegraphics[width=0.4\textwidth,height=0.20\textheight]{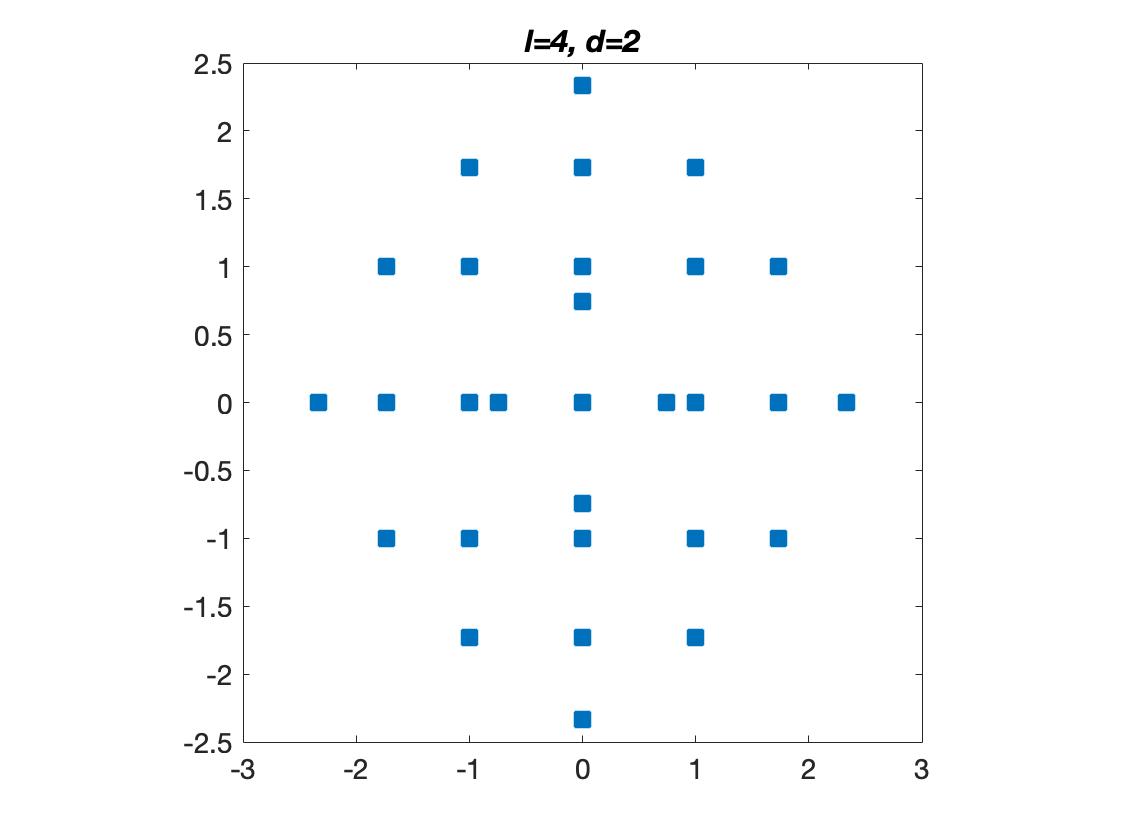}
  \caption{Two-dimensional sparse grid with $l=3$ and $l=4$.}
  \label{fig:sparseGrid2D}
\end{figure}

\begin{figure}[htbp!]
 \centering
 \subfloat[Mean]{%
  \includegraphics[width=0.46\textwidth,height=0.27\textheight]{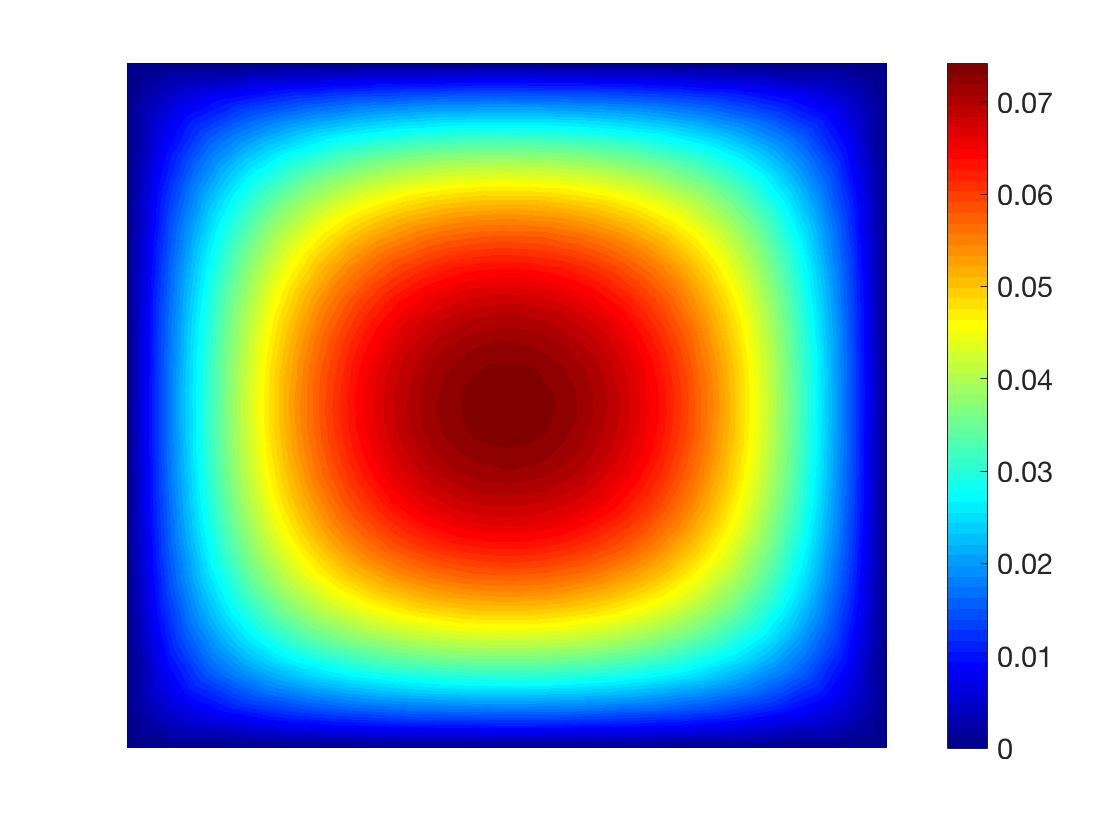}}
  \subfloat[Standard deviation]{%
  \includegraphics[width=0.46\textwidth,height=0.27\textheight]{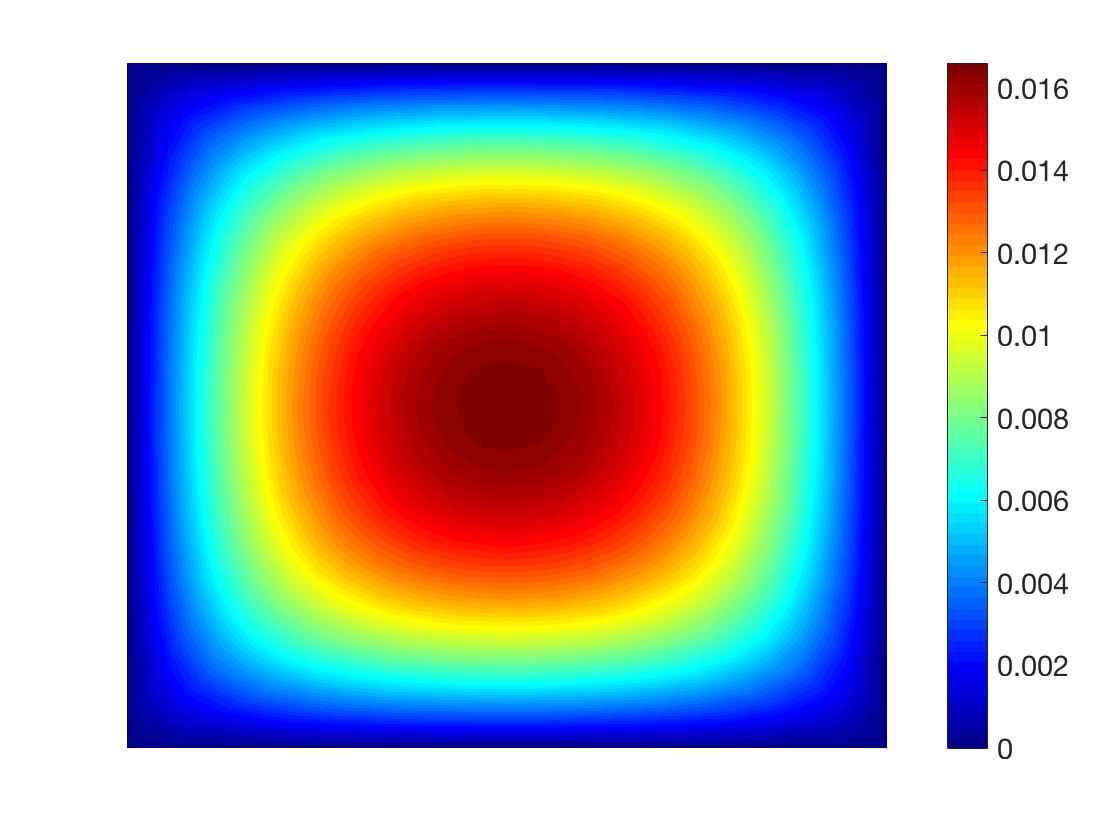}}
  
  \subfloat[${\mathrm{\hat{\bf{u}}}_4}$ \label{subfig:nisp4}]{
  \includegraphics[width=0.46\textwidth,height=0.27\textheight]{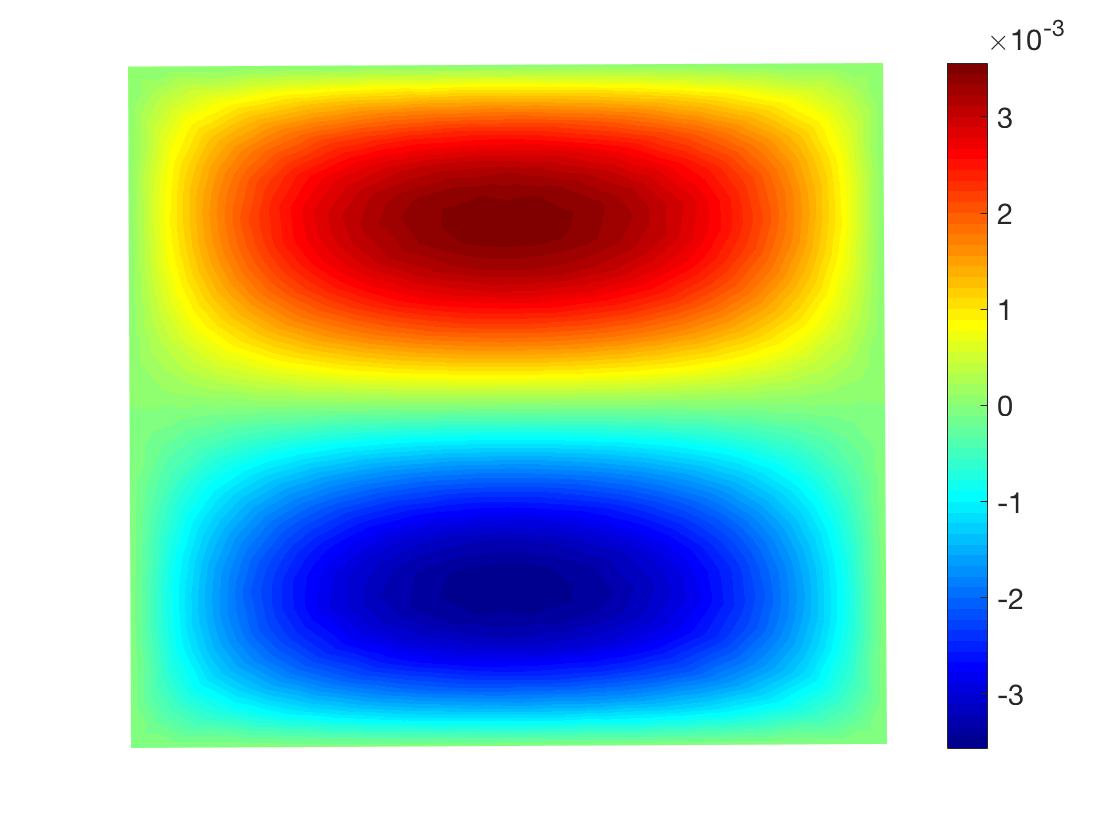}}
  \subfloat[${\mathrm{\hat{\bf{u}}}_5}$ \label{subfig:nisp5}]{
  \includegraphics[width=0.46\textwidth,height=0.27\textheight]{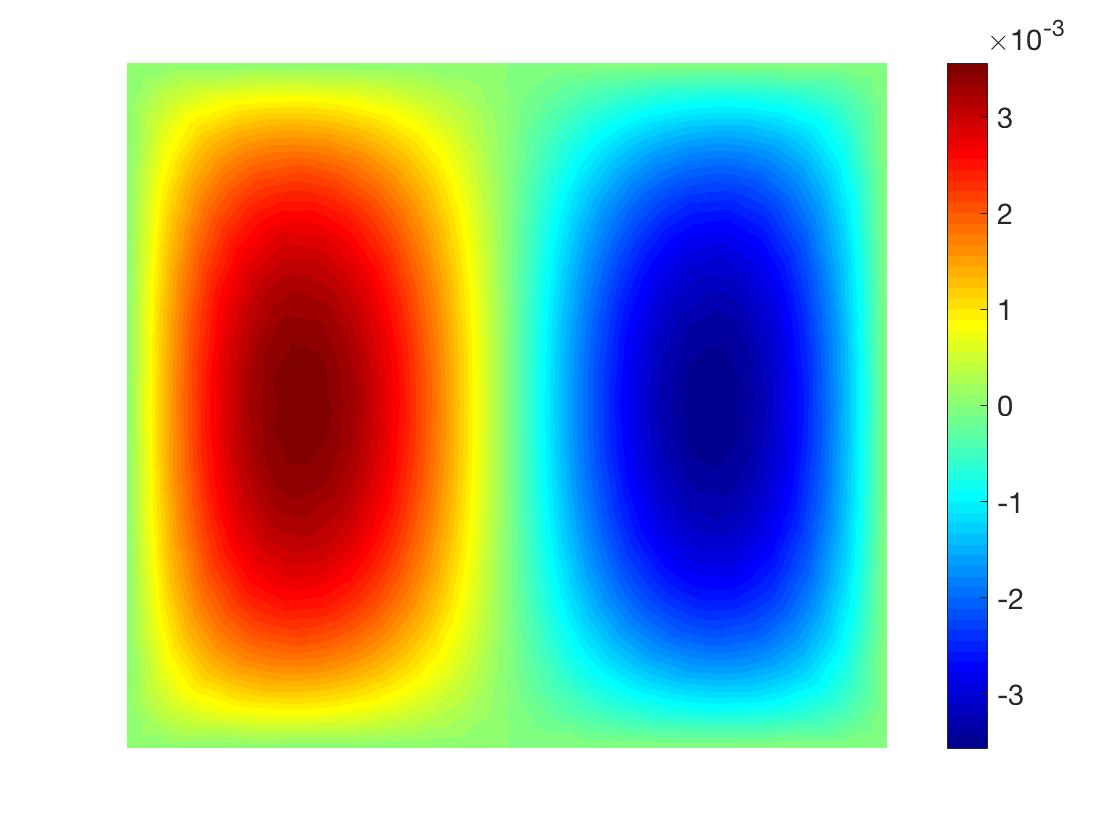}}

  \subfloat[${\mathrm{\hat{\bf{u}}}_{9}}$ \label{subfig:nisp9_2}]{
  \includegraphics[width=0.46\textwidth,height=0.27\textheight]{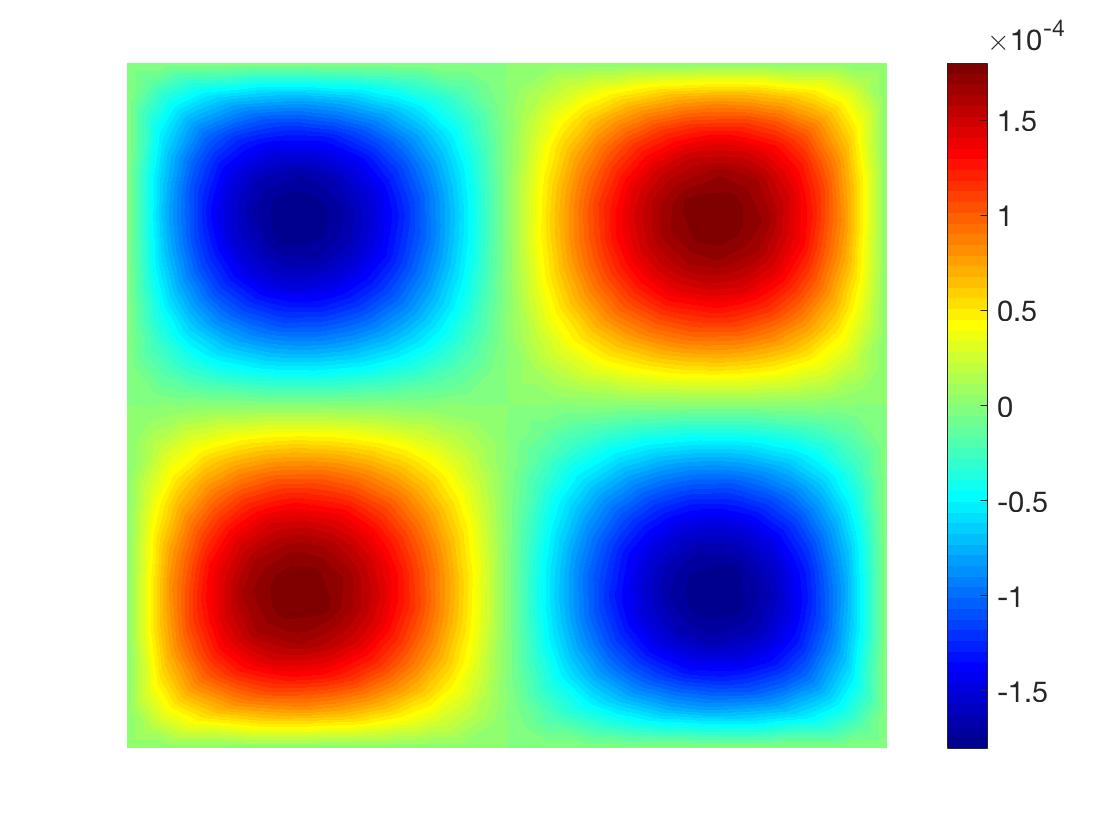}}
  \subfloat[${\mathrm{\hat{\bf{u}}}_{10}}$ \label{subfig:nisp10_2}]{
  \includegraphics[width=0.46\textwidth,height=0.27\textheight]{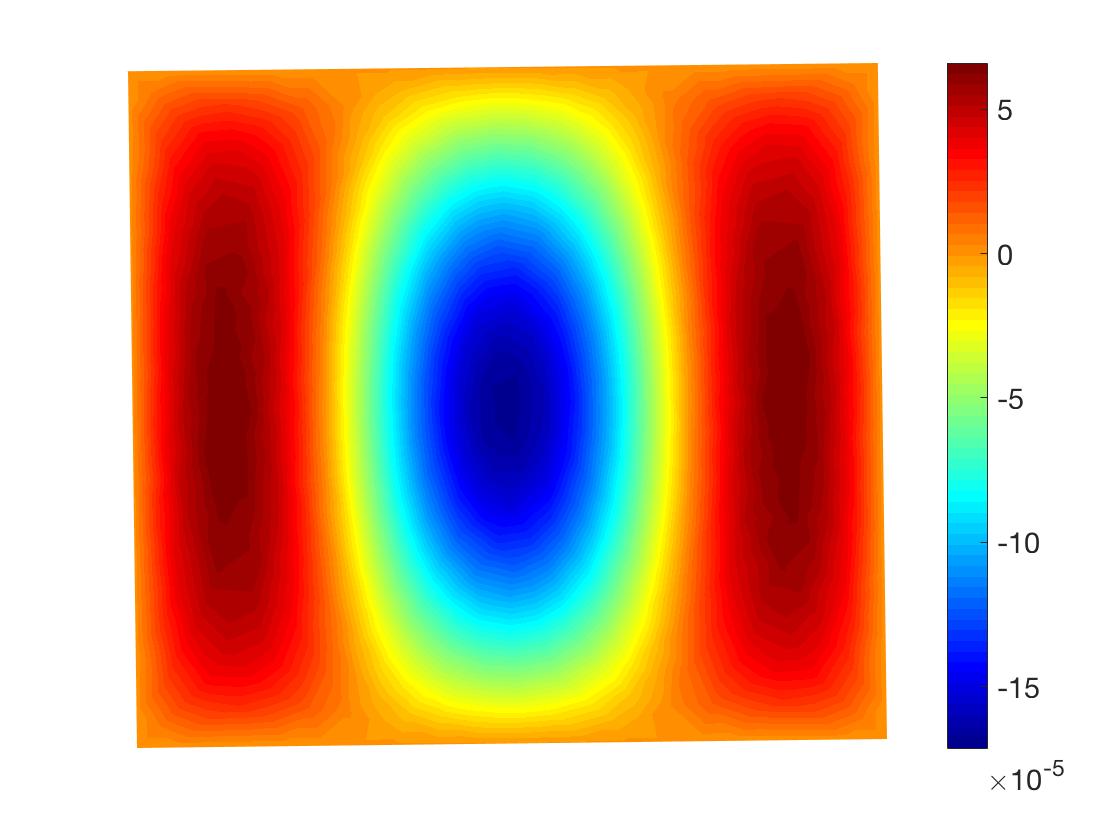}}
  
  \caption{The mean, standard deviation, and a few selected PCE coefficients of the solution field using NISP.}
  \label{fig:nisp_sol}
\end{figure}

The sparse grid in two dimensions for $l=2$ and $l=3$ are shown in \Cref{fig:sparseGrid2D}. Using Smolyak sparse grid quadrature with $l=3$ and $d=2$, deterministic FEM sample evaluations are performed using FEniCS package for the stochastic PDE defined by \Cref{eq:sspe}.
The following numerical parameters $b=1$, $\sigma=0.3$, $L=3$, $p_u=3$ and $f=1$ are used. The procedure outlined in~\Cref{alg:nisp_algo} is employed to obtained PCE coefficients of the solution process.
A unit square domain, discretized using $600$ nodes and $1200$ elements 
is used. The mean, standard deviation, and a few selected PCE coefficients of the solution process are plotted in \Cref{fig:nisp_sol}.

Note that the mean and standard deviation of the solution fields using intrusive (\autoref{fig:intrusiv_sol}) and NISP (\autoref{fig:nisp_sol}) has similar trends. Moreover, the contribution of the higher-order PCE coefficient to the solution process decreases in both intrusive and NISP cases. Among these PCE coefficients, the first order coefficients contain Gaussian contributions and the higher-order coefficients contain the non-Gaussian effects. 

\alglanguage{pseudocode}
\begin{algorithm}
\caption{: \textcolor{blue}{NISP Procedure}}
\label{alg:nisp_algo}
\begin{algorithmic}[1]
\State {\bf{Input} :}  mesh data : $points, edges$ and $triangels$
\For{$ {i = 4, 5, 6, ..., l+1   } $}
\For{$ {j = 2, 3, 4, ..., L}  $}
\State [$nodes$, $weights$] = sparseGrid('Gauss-Hermite', $j$, $i$)
\For{$ {k = 1, ..., length(nodes)}  $}
\State $ck = nodes(k)$;
\State call $b$ = AssembleVector($points$, $edges$, $triangels$, $ck$, 'PDE')
\State call $A$ = AssembleMatrix($points$, $edges$, $triangels$, $ck$, 'PDE')
\State $solve(A u_{k} = b) ==> \{ u \}$
\EndFor
\EndFor
\EndFor
\For{$ {np = 1, ..., P_u}  $}
\State $NR_{np}$ = GaussHermiteQuadrature($weights, u, \Psi_{np}$)
\State $DR_{np}$ = GaussHermiteQuadrature($weights, {\Psi_{np}^2}$)
\EndFor
\State {\bf{Output} :}  ${\mathrm{\hat{\bf{u}}}_{np}} = \frac{NR_{np}}{DR_{np}}$
\end{algorithmic}
\end{algorithm}



\section{Conclusion}\label{sec:conclusion}
In summary, this article demonstrate how to numerically solve stochastic PDEs using FEniCS---a general puprpose deterministic FEM package and UQTk---a collection of libraries and tools for the UQ. The focus is given on the implementation aspects of of intrusive SSFEM in order to reduce the computational complexity arises in the cases of intrusive SSFEM. 

We note that the non-intrusive approach is favorable from an implementational perspective because one can directly employ any existing deterministic solver as a black box to simulate the required samples. On the other hand, the intrusive approach demands additional coding efforts. However, as demonstrated in this article, the stochastic assembly procedure employed for intrusive SSFEM can utilize the readily available deterministic finite element assembly routines (such as FEniCS), which can substantially reduce the coding efforts.

Although we can readily accommodate the increased number of samples due to a large number of random variables in the non-intrusive approach, there is a substantial increase in the number of sample evaluations for the non-Gaussian input making it computationally costly compared to the intrusive approach~\cite{desai2019scalable}. However, the memory required to assemble and solve the intrusive system increases as we increase the number of random variables. Therefore, for a computer with a fixed random access memory (RAM), there is an upper limit to the size of the intrusive system we can accommodate, restricting the application of intrusive SSFEM. 

Nonetheless, when we can handle the increasing intrusive system size by distributing it among multiple nodes and employing an efficient parallel solver, we can solve the intrusive system for solution coefficients much faster than the non-intrusive approach for the same level of accuracy~\cite{desai2019scalable}.
For these reasons, Developing scalable solvers to tackle stochastic PDEs using SSFEM is an active area of research, as evidenced by many articles published in the last couple of decades~\cite{subber2012PhDTh,elman2007solving,mandel2008multispace,ghosh2009feti,stavroulakis2017gpu,pellissetti2000iterative,sousedik2014hierarchical,desai2017scalable,desai2022domain}.

\appendix
\section{Code Snippets}\label{sec:sinppets}

\begin{lstlisting}[caption=L-dimensional polynomials using 1-dimensional polynomials, label=list:multiPoly, style=customm]
% Multidimensional Hermite polynomials (Psi)
% using 1D Hermite polynomials (psi)
%      x  : column vector of evaluation points
%   nord  : order of expansion (p)
%   ndim  : number of dimensions (L)
% nPCTerms: number of PC expansion terms (P)

% Compute non-dimentionalized 1-dimensional Hermite polynomials
% 1-st Order Hermite PC: it is always fixed to psi(1) = 1
psi(:, 1) = 1;
psiD(:, 1) = 1;
% 2nd order Hermite PC: it's always fixed to psi(2) = x
if (nord > 0)
    psi(:, 2) = x;
    psiD(:, 2) = x;
    % 3rd and more are solved by recursive formula
    for i = 3 : nord+1
        psiD(:, i) = x(:) .* psiD(:, i-1) - (i-2) * psiD(:, i-2)
        psi(:, i) = psiD(:, i)/sqrt(factorial(i-1))
    end
end

% Compute non-dimentionalized L-dimensional Hermite polynomials
Psi = ones(1, nPCTerms);
for i = 1 : ndim
    Psi = Psi .* psi(i, multiIndex(1 : nPCTerms, i) + 1);
end
\end{lstlisting}

\

\begin{lstlisting}[caption=Matlab code to compute moments of ND-polynomials using 1D-polynomials, label=list:cijk,  style=customm]
%% Code to calculate moments of N-Dimensional Hermite PC basis.
% Cijk =  <Psi_i Psi_j Psi_k>  : Tripple product require for SSFEM
% input:
%       ndim = stochastic dimension/number of KLE terms (L)
%       nord = order of PCE (p)
%
% output:
%       cijk: <Psi_i Psi_j Psi_k> Tripple product
%       ijk: Non-zero tripple product indices

pcdata1D.ndim = ndim;                  % need to provide this
pcdata1D.nord = nord;                  % need to provide this
pcdata1D.nclp = 2 * nord + 1;          % done as 2n+1 quad points

% Computing 1D quadrature points(x) and weights(w)
[x, w] = quadXW_1D(pcdata1D.nclp);     % refer to UQToolkit
pcdata1D.x = x;
pcdata1D.w = w;

% Non-dimentionalized, 1-dimensional Hermite polynomials
% 1-dimensional Hermite polynomials at quad-points
psi = zeros(nclp, nord + 1);
psiD = zeros(nclp, nord + 1);
% 1-st order Hermite PC: it is always fixed to psi(1) = 1
psi(:, 1) = 1;
psiD(:, 1) = 1;
% 2nd order Hermite PC: it's always fixed to psi(2) = x
if (nord > 0)
    psi(:, 2) = x;
    psiD(:, 2) = x;
    % 3rd and more order Hermite PC's are solved by recursive formula
    for i = 3 : nord+1
        psiD(:, i) = x(:) .* psiD(:, i-1) - (i-2) * psiD(:, i-2)
        psi(:, i) = psiD(:, i)/sqrt(factorial(i-1))
    end
end

% multi-indices and the number of terms in PCE: refer to UQToolkit
[pcdata1D.multiIndex, pcdata1D.nPCTerms] = multiIndex(nord, ndim);
% 1D Moments of Hermite PC basis
apow=zeros(nord + 1, nord + 1, nord + 1);
for k = 1 : nord + 1
    for j = 1 : nord + 1
        for i = 1 : nord + 1
            sum = 0;
            for m = 1 : nclp
                sum = sum + psi(m,i)*psi(m,j)*psi(m,k)*w(m);
            end
            apow(i,j,k) = sum;
        end
    end
end

% KLE data: Select desired number of KLE-order and KLE-dim
KLEord = 2;
KLEdim = ndim;
PCEord = 3;
nPCEin = factorial(KLEord+KLEdim)./(factorial(KLEord).*factorial(KLEdim))
nPCTerms = factorial(PCEord+KLEdim)./(factorial(PCEord).*factorial(KLEdim))
% ND Moments using 1D moments
tol = 1e-8;
indexi = 1;
for i = 1 : nPCTerms
    for j = 1 : nPCTerms
        aprod = 1;
        for m = 1 : ndim
            l1 = multiIndex(i, m) + 1;
            l2 = multiIndex(j, m) + 1;
            l3 = multiIndex(k, m) + 1;
            aprod=aprod*apow(l1, l2, l3);
        end
        temp1 = aprod;
        if (temp1 > tol)
            indexi;
            ijk = [ijk, [k;i;j]];
            cijk = [cijk, temp1];
        end
        indexi = indexi+1;
    end
end
cijk = roundn(cijk, -6);
nonZeroIndex = ijk;
\end{lstlisting}

\

\begin{lstlisting}[caption=FEniCS based python code snippet for deterministic FEM, label=list:det_fem, style=customp]
## Code snippet to solve Poisson Equation with Dirichlet BC
from dolfin import *
parameters['reorder_dofs_serial'] = False

# import (external using GMSH) mesh
mesh = Mesh("unitSquare.xml")
V = FunctionSpace(mesh, "Lagrange", 1)

# Define boundary conditions
def boundary(x):
    return x[0] < DOLFIN_EPS or x[0] > 1.0 - DOLFIN_EPS or
           x[1] < DOLFIN_EPS or x[1] > 1.0 - DOLFIN_EPS
u0 = Constant(0.0)
bc = DirichletBC(V, u0, boundary)

# Define variational problem
u = TrialFunction(V)
v = TestFunction(V)
f = Constant(1.0)
cd = Constant(1.05)
a = Form(cd*inner(grad(u), grad(v))*dx)
L = Form(f*v*dx + g*v*ds)

# Compute solution
u = Function(V)
solve(a == L, u, bc)

# Save solution to VTK format
file = File(poissonOriginal.pvd)
file << u
\end{lstlisting}

\

\begin{lstlisting}[caption=FEniCS based intrusive SSFEM assembly code snippet, label=list:sto_femAssembly, style=customp]
## Code for intrusive system matrix assembly using FEniCS-based,
## element-level deterministic assembly routines
from dolfin import *
parameters['reorder_dofs_serial'] = False

# import (external using GMSH) mesh
mesh = Mesh("unitSquare.xml")
# Initialize the connectivity between facets and cells
tdim = mesh.topology().dim()
mesh.init(tdim-1, tdim)

# Define variational problem
V = FunctionSpace(mesh, "Lagrange", 1)
u = TrialFunction(V)
v = TestFunction(V)
f = Constant(1.0)

## Load pre-calculated ijk and Cijk files
cijkMat = sio.loadmat("cijk.mat")
ijk = cijkMat['ijk']
cijk = cijkMat['cijk']

# Input/Output PCE data
PCE_A = 10     #p=2, L=3
PCE_u = 20     #p=3, L=3

indexi = 0
N = V.dim()              ## deterministic system matrix size
Ns = (PCE_u*N)           ## stochastic system matrix size
As = np.zeros([Ns,Ns])   ## initialize intrusive system matrix
for inputPCE_index in range(PCE_A):
    ## Call stochastic variational form for each PCE input
    a, L = stoVariationalForm(inputPCE_index, u, v, g, f)
    ## Call deteministic assembly for each PCE input
    A, b = detFemAssembly(a, L, inputPCE_index)

    for i in range(PCE_u):
        for j in range(PCE_u):
            if inputPCE_index==ijk[indexi,0] and i==ijk[indexi,1] and j==ijk[indexi,2]:
                As[i*N:(i+1)*N, j*N:(j+1)*N] += cijk[indexi]*A
                As_old[i*N:(i+1)*N, j*N:(j+1)*N] += cijk[indexi]*A_old
                indexi = indexi+1

## Save assembly matrices in ".mat" or ".txt" format
sio.savemat('As.mat', {'As':As, 'b':b})
# np.savetxt('As,dat', As), np.savetxt('bs,dat', b)
\end{lstlisting}

\begin{lstlisting}[caption=Python code snippet for FEniCS based (element-level) FEM assembly , label=list:det_femAssembly, style=customp]
## Function to perform FEniCS-based, element-level deterministic assembly
## Procedure is employed for stochastic system matrix assembly

def detFemAssembly(a, L, inputPCE_index):
    # Dummy problem to define LAYOUT of global problem
    A_g = assemble( Constant(0.)*inner(grad(u), grad(v))*dx )
    f_g = assemble( Constant(0.)*v*dx )

    # Get dofmap to construct cell-to-dof connectivity
    dofmap = V.dofmap()

    # Perform assembly
    for cell in cells(mesh):
        dof_idx = dofmap.cell_dofs(cell.index())
        # Assemble local rhs and lhs system
        a_local  = assemble_local(a, cell)
        L_local  = assemble_local(L, cell)
        # Assemble global rhs and lhs system
        A_g.add_local(a_local,dof_idx, dof_idx)
        f_g.add_local(L_local,dof_idx)

    # Finalize assembling
    A_g.apply("add"), f_g.apply("add")

    # Define Dirichlet boundary condition
    def boundary(x):
    return x[0] < DOLFIN_EPS or x[0] > 1.0 - DOLFIN_EPS or
        x[1] < DOLFIN_EPS or x[1] > 1.0 - DOLFIN_EPS
    u0 = Constant(0.0)
    bc = DirichletBC(V, u0, boundary)

    # Apply boundary condition
        if inputPCE_index==0:
        bc.apply(A_g,f_g)

    # Get values
    A = A_g.array()
    b = f_g.getlocal()

    return A, b
\end{lstlisting}

\

\begin{lstlisting}[caption=Python code snippet for FEniCS based stochastic variational form of Poisson equation, label=list:sto_variation, style=customp]
## FEniCS-based variational form for stochastic diffusion equation
## Procedure is employed for stochastic system matrix assembly
def stoVariationalForm(inputPCE_index, u, v, f):
    cd = MyExpression(inputPCE_index, degree=0)
    a = Form(cd*inner(grad(u), grad(v))*dx)
    L = Form(f*v*dx)
    return a, L
\end{lstlisting}

\

\begin{lstlisting}[caption=Python code snippet for PC expansion of lognormal stochastic process, label=list:sto_cd, style=customp]
# Function class to get PCE coefficients for c_d = exp(g)
# c_d is defined as a lognormal stochastic process
# inputs : takes four inputs [index, ndim, sIndex, mIndex]
# output : one c_d coefficient value for each inputPCE_index


class MyExpression(Expression):
    def __init__(self, inputPCE_index, **kwargs):
        self.index = inputPCE_index[0]     # PCE_index
        self.ndim = ndim            # Number of RVs (need to provide)

        sIndexPath = "sortIndex.dat"  # load precalculated sort-index
        sIndex = numpy.genfromtxt(sIndexPath)
        mIndexPath = "multiIndex.dat"  # load precalculated multi-index
        mIndex = numpy.genfromtxt(mIndexPath)
        self.sIndex = params[2].astype(int)  # sort-index
        self.mIndex = params[3].astype(int)  # multi-index

    # Evaluate c_d at x
    def eval(self, values, x):
        # Precalculated omegas and lambdas for exponential covariance
        # kernal with bx=by=b=1, sigma=0.3 on a unit square domain
        # multipliers = sqrt(lambda)/sqrt(a-(sin(2*w*a)/(2*w)))
        multipliers = [0.92184, 0.49248, 0.29374, 0.20437, 0.15576]
        omegas = [1.30654, 3.67319, 6.58462, 9.63168, 12.72324]

        # Mean and standard deviation of underlying Gaussian process
        meang = 1.05  # For mu=0, sigma=0.3
        sigma = 0.3

        # Trunctaed PCE of lognormal process: Automated to n-RVs
        # KLE: obtain KLE coefficients for the input
        g = []
        for i in range(self.ndim):
            Xindex = self.sIndex[i, 0]
            if (Xindex % 2) == 0:  # even
                gg1 = multipliers[Xindex-1] * (sin(omegas[Xindex-1]*(x[0]-0.5)))
            else:  # odd
                gg1 = multipliers[Xindex-1] * (cos(omegas[Xindex-1]*(x[0]-0.5)))

            Yindex = self.sIndex[i, 1]
            if (Yindex % 2) == 0:  # even
                gg2 = multipliers[Yindex-1] * (sin(omegas[Yindex-1]*(x[1]-0.5)))
            else:  # odd
                gg2 = multipliers[Yindex-1] * (cos(omegas[Yindex-1]*(x[1]-0.5)))

            g.append(sigma*gg1*gg2)

        # PCE: obtain PCE coefficients for the (lognormal) input
        Y = 1.0
        i = self.index
        for j in range(self.ndim):
            idx = mIndex[i, j]
            if idx == 0:
                yy = 1.0
            elif idx == 1:
                yy = g[j]
            else:
                nfactorial = factorial(idx)  # Use factorial function
                yy = (g[j]**(idx))/numpy.sqrt(nfactorial)

            Y = Y*yy

        values[0] = meang * Y
        # PCE coefficient for the input term
\end{lstlisting}



\bibliographystyle{unsrt}  
\bibliography{references}


\end{document}